\documentclass{CSML}
\pdfoutput=1

\usepackage{lastpage}
\newcommand{\dOi}{13(4:13)2017}
\lmcsheading{}{1--31}{}{}%
{Dec.~24, 2015}{Nov.~21, 2017}{}

\usepackage{hyperref}
\hypersetup{hidelinks}

\usepackage[]{amsmath, amssymb}

\begin{document}

\title[Computing Solution Operators of PDEs]{Computing Solution Operators of \\ Boundary-value Problems
for \\ Some Linear Hyperbolic Systems of PDEs}

\author{Svetlana Selivanova}
\address{S.L. Sobolev Institute of Mathematics
\\Siberian Branch of the
Russian Academy of Sciences; Novosibirsk State University\\
    Novosibirsk, Russia} 
\email{s\_seliv@math.nsc.ru}
 \author{Victor Selivanov}
\address{A.P. Ershov Institute of
Informatics Systems\\Siberian Branch of the Russian Academy of
Sciences; Novosibirsk State University\\
    Novosibirsk, Russia}
    \email{vseliv@iis.nsk.su}


\begin{abstract}
We discuss possibilities of application of Numerical Analysis methods
 to proving computability, in the sense of the TTE approach,
 of solution operators of boundary-value problems for systems of PDEs.
 We prove computability of the solution operator for a
symmetric hyperbolic system with computable real coefficients and
dissipative boundary conditions, and of the Cauchy problem for the
same system (we also prove computable dependence on the
coefficients) in a cube $Q\subseteq\mathbb R^m$. Such systems
describe a wide variety of physical processes (e.g.\ elasticity,
acoustics, Maxwell equations). Moreover, many boundary-value
problems for the wave equation also can be reduced to this case,
thus we partially answer a question raised in \cite{wz02}. Compared
with most of other existing methods of proving computability for
PDEs, this method does not require existence of explicit  solution
formulas and is thus applicable to a broader class of (systems of)
equations.
 \end{abstract}

 \maketitle

\section{Introduction}\label{in}

We consider  boundary-value problems for systems of PDEs of the
form
\begin{equation} \label{sist}
\begin{cases} L
{\bf u}(y)=f(y)\in C^p(\Omega,\mathbb R^n),\quad y\in\Omega\subset\mathbb R^{k}\\
{\mathcal L}{\bf u}(y)|_{\Gamma}=\varphi(y\mid_{\Gamma})\in
C^q(\Gamma,\mathbb R^n),\ \Gamma\subseteq \partial\Omega,
\end{cases}
\end{equation}
where $L$ and ${\mathcal L}$ are differential operators (the
differential order of ${\mathcal L}$ is less than the one of $L$),
$\Gamma$ is a part of the boundary $
\partial\Omega$ of some area $\Omega$. In particular, if
$\Gamma=\{t=0\}$ and $t$ is among  the variables
$y_1,y_2,\ldots,y_k$, then~\eqref{sist} is a Cauchy (or
initial-value) problem.
Assuming existence and uniqueness of the solution ${\bf u}$ in
$\Omega$, we study computability properties of the solution operator
$R:(L,\mathcal{L},f,\varphi)\mapsto{\bf u}$.  Note that in~
\eqref{sist} the number $k$ of ``space'' variables
$y_1,y_2,\ldots,y_k$ is not necessarily equal  to the number $n$ of
the unknown functions $u_1,u_2,\ldots,u_n$, e.g.\ for the linear
elasticity equations~\eqref{elast} we have $n=9,k=4$.

Computability will be understood in the sense of Weihrauch's TTE
approach \cite{wei}. Recently the following main achievements in
the study of computability properties of PDEs were made.
Computability of solution operators of initial-value problems for
the wave equation \cite{wz02}, Korteveg de Vries equation
\cite{gzz01,wz05}, linear and nonlinear Schr\"odinger equations
\cite{wz06} was established; also  computability of fundamental
solutions of PDEs with constant coefficients
$Pu=\sum\limits_{|\alpha|\leq M}c_{\alpha}D^{\alpha}u=f$ was
proved in \cite{wz06-2}. Most of the methods of the mentioned papers are
based on a close examination of explicit solution formulas and the
Fourier transformation method, except for the paper \cite{wz05}
where a method based on fixed point iterations is introduced. In
these papers, the initial data and solutions are mainly assumed to
belong to some Sobolev classes of generalized functions.

As  is well-known, explicit solution formulas for boundary-value
problems (even for the Cauchy initial-value problems) exist rarely.
Even for the simplest example of  the wave equation the
computability of the solution operator for boundary-value problem
was formulated in \cite{wz02} as an open question, and we have not
seen any paper where this question would be answered. Results of our
paper provide, in particular, a positive answer to this question
for the case of computable real coefficients and dissipative
boundary conditions, for classes  of continuously differentiable
functions with uniformly bounded derivatives.

In \cite{ss09} we propounded an approach to study the
computability of PDEs based on finite-dimensional approximations
(difference schemes widely used in numerical analysis) and established
computability, in the sense of the TTE approach, of the solution
operator $\varphi\mapsto {\bf u}$ of the Cauchy problem for a
symmetric hyperbolic system, with a zero right-hand part, of the
form
\begin{equation} \label{sist_1}
\begin{cases} A\frac{\partial
{\bf u}}{\partial t}+\sum\limits_{i=1}^mB_i\frac{\partial {\bf u}}
{\partial x_i}=0,\ t\geq 0,\\
{\bf u}|_{t=0}=\varphi(x_1,\ldots,x_m).
\end{cases}
\end{equation}
Here $A=A^\ast>0$ and $B_i=B_i^\ast$ are constant symmetric
computable $n\times n$-matrices, $t\geq0$, $x=(x_1,\ldots,x_m)\in
Q=[0,1]^m$, $\varphi:Q\rightarrow{\mathbb R}^n$ and ${\bf
u}:Q\times[0,+\infty)\rightharpoonup{\mathbb R}^n$ is a partial
function acting on the domain $H$ of existence and uniqueness of
the Cauchy problem~\eqref{sist_1}. In \cite{ss09} the
computability of the domain $H$ (which is a convex polyhedron
depending only on $A,B_i$) was also proved. The operator $R$
mapping a $C^{p+1}$ function $\varphi$ to the unique $C^p$
solution ($p\geq 2$) is computable, if the norms of the first and
second partial derivatives of $\varphi$ are uniformly bounded.

Such systems   can be used to describe a wide variety of physical
processes like those considered in the theories of elasticity,
acoustics, electromagnetism etc., see e.g.\ 
\cite{fri,go71,go76,ll86,ll04,kps,gm98}. They were first considered in 1954 by
K.O. Friedrichs  \cite{fri}. He proved the existence theorem
based on finite difference approximations, in contrast with the
Schauder-Cauchy-Kovalevskaya method based on approximations by
analytic functions and a careful study of infinite series. The
notion of a hyperbolic system (applicable also to broader classes
of systems) is due to I.G. Petrovskii \cite{petr}, see also the
very interesting discussion on different notions of hyperbolicity
and their motivations in \cite{fri}.

Recall that a linear first-order differential operator
$E=\sum\limits_{\mu=1}^mA_{\mu}\frac{\partial}{\partial
x_{\mu}}+B,$ where $A_{\mu}$, $B$ are real $n\times n$ matrices,
$\mu=1,2,\ldots,m$, is called {\it hyperbolic} in the sense of
Petrovskii, if there is a
$\xi^0\in\mathbb R^m$ such that, for all $\xi\in\mathbb R^m$, the
matrix pencil
$$\sum\limits_{\mu=1}^m\xi_{\mu}A_{\mu}-\lambda\sum\limits_{\mu=1}^m\xi^0_{\mu}A_{\mu}$$
has
 real eigenvalues $\lambda$.
In particular, if all the matrices $A_{\mu}$, $\mu=1,2,\ldots,m$
are symmetric and one of them is positive-definite, as in~
\eqref{sist_1}, then the operator $E$ is obviously hyperbolic in
this sense.

The Friedrichs' method has turned out to be interesting from the
computational point of view because it yields algorithms for
solving PDEs's in the exact sense of Computable Analysis  which
are based on methods really used in Numerical Analysis.

In this  paper we  prove computability for a broad class of
boundary-value problems for~\eqref{sist_1}, by using the
difference approximations approach stemming from the work
\cite{fri} and developed in \cite{gr,go71,go76,kps} and others.
Many details of our proofs are similar to those of the proof of
the existence theorem for the linear hyperbolic systems in
\cite{go71,go76,fri} but, since we refer to more rigorous approach of
computable analysis we are forced to establish several additional
estimates. Actually, these proofs are based on careful estimates of the difference approximations of the considered differential operators, ideas of which can be also found e.g.\ in  
\cite{fri,strik,gr,gv}.

Our study intensively uses the well-known classical theorem of the theory
of difference schemes stating that the approximation and stability
properties of a difference scheme imply its convergence to the
solution of the correspondent differential equation in a suitable grid
norm uniformly on steps.

The  proofs of this paper rely also on  the well-known fact that the
ordered field of algebraic real numbers  and some extensions of this
field are strongly constructivizable (this is closely related to the
Tarski's quantifier elimination for real closed fields, see e.g.\ 
\cite{ta51,bp06}) which implies  computability of necessary spectral
characteristics of symmetric matrices with algebraic real
coefficients. This makes obvious  computability of all steps in
the iterative process induced by the difference scheme used in this
paper. This trick also leads to an improvement  of the main result in \cite{ss09} to
the result that the solution operator for the Cauchy problem~
\eqref{sist_1} is computable not only on $\varphi$ but also on the
coefficients  $A,B_i$ (under some additional assumptions, see
Theorem~\ref{main2}).

Our proofs here give some additional details compared to those in  \cite{ss09}. They  make use of results in several fields: PDEs,
difference schemes, computable analysis, computable fields. The results and proofs establish a close connection between computable number fields and computable reals and apply this connection to proving computability of solutions of some PDEs. 
Unfortunately, they do not yield reasonable upper bounds for the complexity of
solving the initial and boundary value problems for PDEs because the corresponding results on computable fields are established here with the use of unbounded search algorithms. Search for more feasible algorithms is a natural further step in the study of
computability properties of PDEs.

In Section~\ref{state} we describe the considered problems and
assumptions we need to prove the  computability of solution
operators. Some necessary notions and facts are recalled in
Section~\ref{prel}. Section~\ref{fdim} is devoted to the
construction of a difference operator approximating the
differential problem and its basic properties. In Section~\ref{smain}
we formulate precisely the main results of the paper
and describe the proof schemes, without technical details of the
corresponding estimates. The technical details are proved in
Section~\ref{estimates}. We conclude in Section~\ref{concl} by a
short discussion on more general systems~\eqref{sist}.

\section{Statement of the boundary-value problem and examples}\label{state}

Along with the Cauchy problem~\eqref{sist_1} we now consider the
following boundary-value problem:
\begin{equation}\label{sist_2}
\begin{cases} A\frac{\partial
{\bf u}}{\partial t}+\sum\limits_{i=1}^mB_i\frac{\partial {\bf u}}
{\partial x_i}=f, \\
{\bf u}|_{t=0}=\varphi(x_1,\ldots,x_m), \\
\Phi_i^{(1)}{\bf u}(x_1,\ldots,x_{i-1},0,x_{i+1},\ldots,x_m,t)=0,\\
\Phi_i^{(2)}{\bf u}(x_1,\ldots,x_{i-1},1,x_{i+1},\ldots,x_m,t)=0,\\
i=1,2,\ldots,m,
\end{cases}
\end{equation}
where
\begin{itemize}
\item $A=A^\ast>0$ is positively definite and $B_i=B_i^\ast$  are  fixed computable symmetric
$n\times
n$-matrices;
\item $0\leq t\leq T$ for a computable real $T$;
\item $x=(x_1,\ldots,x_m)\in Q=[0,1]^m$;
\item $\varphi\in C^{p+1}(Q,\mathbb R^n)$,  $f\in
C^{p}(Q\times[0,T],\mathbb R^n)$, $p\geq 2$ (in this paper we let $f=0$ for simplicity);
\item the boundary coefficients $\Phi_i^{(1)}$, $\Phi_i^{(2)}$ are fixed
computable  matrices meeting the
following conditions:

\begin{enumerate}[label=\arabic*)]
\item The number of rows of $\Phi_i^{(1)}$ (respectively,
$\Phi_i^{(2)}$) is  equal to the number of positive (respectively,
negative) eigenvalues of the matrices $A^{-1}B_i$, and the boundary values of ${\bf u}$ are  consistent  with the initial conditions $\varphi$;

\item The boundary conditions are assumed to be  dissipative which
means that
\begin{equation}\label{gran_dissip}
\langle B_i{\bf u},{\bf u}\rangle\leq 0\text{ for }x_i=0,\quad \langle B_i{\bf u},{\bf
u}\rangle \geq 0\text{ for }x_i=1,\quad i=1,2,\ldots,m.
 \end{equation}
\end{enumerate}
\end{itemize}

\begin{rem}Note that the assumptions 1) regarding the dimensions of the matrices
$\Phi_i$ and consistency of the boundary conditions with the initial ones are needed for proving existence
 of a solution ${\bf u}\in C^p(Q\times[0,T], \mathbb R^n)$ of~\eqref{sist_2},
while the assumption~\eqref{gran_dissip} provides uniqueness of
the solution \cite{fri,go71,evans,jo66}. 

Moreover, these assumptions are
needed \cite{go76,gr} for proving stability of the
 difference scheme constructed below in Section~\ref{fdim},
 which is one of the main ingredients in the proof of computability
 results. 
 
 The smoothness assumptions $\varphi\in C^{p+1}(Q,\mathbb R^n)$,  $f\in
C^{p}(Q\times[0,T],\mathbb R^n)$, $p\geq 2$,    are needed to provide  at least $C^2$ smoothness of the solution, which is essential to establish estimates of the convergence constant in what follows.

The consistency conditions have to be found for each particular case of the boundary conditions and matrix coefficients and are usually nontrivial. We don't go into details since the algorithms below do not depend on their concrete expressions. They matter only for proofs of the existence and uniqueness theorems.

\end{rem}

An example of dissipative boundary conditions for the system~
\eqref{sist_2} are conservative boundary conditions, stating that
the energy flow through the boundary is constant:
\begin{equation}\label{energy1}\oint\limits_S \langle[\tau A+
\sum\limits_{i=1}^m\xi_i B_i]{\bf u},\textbf{u}\rangle dS=
\int\limits_{Q\times[0,T]}2\langle f,{\bf u}\rangle d\Omega,\end{equation} where
$(\tau,\xi_1,\ldots,\xi_m)$ is the extrinsic normal vector for the
surface $S=\partial (Q\times[0,T])$.

E.g.\ for the linear elasticity equations (which constitute a
symmetric hyperbolic system with $9\times 9$-matrices)
\begin{equation}\label{elast}
 \begin{cases}
 \frac{1}{2\mu}\frac{\partial \sigma _{ij}}{ \partial t}-\frac{\lambda}
{2\mu(3\lambda+2\mu)}\delta _{ij}\frac{\partial(\sigma
_{11}+\sigma_ {22}+ \sigma _{33})} {\partial t}-
 \frac{1}{2}(\frac{ \partial u_i } {\partial x_j}+
\frac{\partial u_j} {\partial x_i}) =0, \quad i,j=1,2,3,\vspace{6.pt}\\
  \rho\frac{\partial u_i}{
\partial t}- \frac{\partial\sigma _{ij} } { \partial x_j}=0,
\quad i=1,2,3,
\end{cases}\end{equation}
in particular, the following boundary equations are conservative
 $$\begin{cases}
x_1=0,\ x_2=1:\ \sigma_{11}=\sigma_{12}=\sigma_{13}=0,\\
y_1=0,\ y_2=1:\ \sigma_{12}=\sigma_{22}=\sigma_{23}=0,\\
z_1=0,\ z_2=1:\ \sigma_{13}=\sigma_{23}=\sigma_{33}=0,
\end{cases}$$
which means that the tangent stresses at any boundary are zero.
Indeed, the energy conservation law~\eqref{energy1} takes the form
 $$\oint\limits_S [-2\xi (u_1\sigma _{11}+
u_2\sigma _{12}+u_3\sigma _{13})-2\eta(u_1\sigma _{12}+ u_2\sigma
_{22}+u_3\sigma _{23})-2\zeta(u_1\sigma _{13}+ u_2\sigma _{23}+u_3\sigma
_{33})]dS=0.$$ Here $u_i$ are the velocities, $\sigma_{ij}$ is the
stresses tensor (a symmetric $3\times 3$--matrix with  6 independent
variables), $\rho$ is the density and $\lambda, \mu$ are the Lame
coefficients.

Interestingly, the boundary-value problem for the wave equation
\begin{equation}\label{wave}
\large\begin{cases} p_{tt}-c_0^2(p_{xx}+p_{yy}+p_{zz})=0,\ \
(x,y,z)\in
\Omega=[0,1]^3,\\
p|_{t=0}=\varphi(x,y,z),\\
p_t|_{t=0}=\psi(x,y,z),\\
p|_{\partial\Omega}=0.\end{cases}
\end{equation}
can be reduced, in several  ways, to a symmetric hyperbolic system
(see e.g.\ \cite{gor79,evans,jo66}), in particular to the
three-dimensional acoustics equations
{\begin{equation}\label{sp_ak}\begin{cases}
\rho_0\frac{\partial u}{\partial t} + \frac{\partial p}{\partial x}=0,\\
\rho_0\frac{\partial v}{\partial t} + \frac{\partial p}{\partial y}=0,\\
\rho_0\frac{\partial w}{\partial t} + \frac{\partial p}{\partial z}=0\\
\frac{\partial p}{\partial t}+\rho_0c_0^2(\frac{\partial
u}{\partial x}+
\frac{\partial v}{\partial y}+\frac{\partial w}{\partial z})=0,\\
p|_{t=0}=\varphi(x,y,z),\\
u|_{t=0}=-\frac{1}{\rho_0c_0^2}\int\limits_0^x\psi(\xi,y,z)d\xi,\\
v|_{t=0}=0,\\
w|_{t=0}=0,\\
p|_{\partial\Omega}=0,
\end{cases}\end{equation}

where $u,v,w$ are the velocities, $p$ is the pressure,
$\rho_0$ is the density and $c_0$ is the speed constant.

Obviously, such a reduction can be done effectively, since
integration is a computable operation. Thus the methods of proving
computability for symmetric hyperbolic systems can be also applied
to prove computability for the wave equation. This gives a partial
answer  to an open question raised in
\cite{wz02}: a  boundary-value problem for the wave equation is
computable (in classes of functions  with uniformly bounded
derivatives, see the exact formulation below) provided that it is
dissipative (i.e., the corresponding boundary-value problem for a
symmetric hyperbolic system to which the wave equation is reduced,
is dissipative) and $c_0$ is a computable real.

\bigskip

We prove  computability of the solution
operator $\varphi\mapsto{\bf u}$ of the boundary-value problem~
\eqref{sist_2} under the following additional assumptions/modifications.

\begin{itemize}
\item First of all, note, that the cube $Q=[0,1]^m$ can easily
be replaced by a computable parallelepiped
 $$[x_1^{(1)},x_2^{(1)}]\times[x_1^{(2)},x_2^{(2)}]\times\ldots\times[x_1^{(m)},x_2^{(m)}],$$
and, in place of $t\geq 0$, we can assume that $t\geq t_0$, where
$t_0$ is a computable real.

\item The first and second partial derivatives of the initial
 function $\varphi$ are bounded by a uniform constant.

\item The considered spaces $C^p$ are equipped with the
$sup$-norm on $Q$ or the $sL_2$-norm on $Q\times[0,T]$, which is an
$L_2$-norm over $Q$ and a $sup$-norm over $[0,T]$, see the more
precise definitions in the next section.

\item For the Cauchy problem~\eqref{sist_1}, we also prove
computability of the solution on the matrices $A,B_i$ assuming
them to belong to the set of symmetric matrices with $A>0$,
the norms
 $$||A||_2=\lambda_{max}(A),
||A^{(-1)}||_2=\frac{1}{\lambda_{min}(A)},\; ||B_1||_2,
||B_2||_2,\ldots, ||B_m||_2$$
  to be bounded by a uniform constant,
the matrix  pencils $\lambda A- B_i$ to have no zero eigenvalues, and
to have the cardinalities of spectra (as well as the cardinality of spectrum of the matrix $A$) given as inputs. Here
$\lambda_{max}(A)$ and $\lambda_{min}(A)$ are the maximal and
minimal eigenvalues of $A$, respectively. This result improves the
main result in \cite{ss09}.
\end{itemize}

\section{Preliminaries}\label{prel}

In this section we briefly  summarize some relevant notions and
facts. In Subsection~\ref{compcount}
we briefly recall some basic notions and facts about constructive
structures, with an emphasis on computable fields. Subsection~\ref{compu}
contains some relevant information on computable
metric spaces. Concrete metric spaces relevant to this paper are
carefully described in Subsection~\ref{spaces}.

\subsection{Computability on countable structures}\label{compcount}

We briefly  recall some relevant notions and facts from computable
model theory.

Recall that a {\em numbering} of a set $B$ is a surjection $\beta$
from $\mathbb{N}$ onto $B$. For numberings $\beta,\gamma$ of $B$,
$\beta$ is {\em reducible} to $\gamma$ (in symbols
$\beta\leq\gamma$) iff $\beta=\gamma\circ f$ for some computable
function $f$ on $\mathbb{N}$, and $\beta$ is {\em equivalent} to
$\gamma$ (in symbols $\beta\equiv\gamma$) iff $\beta\leq\gamma$ and
$\gamma\leq\beta$.  These notions (introduced by A.N. Kolmogorov)
enable to transfer the computability theory over $\mathbb{N}$ to
computability theory over  many other countable structures. The
notions apply to arbitrary functions $\beta,\gamma:\mathbb{N}\to B$
(not only surjections).  Natural extensions of these notions to
partial numberings (i.e., functions defined on a subset of
$\mathbb{N}$) are also of use. In this case, $\beta\leq\gamma$ means
the existence of a computable partial function $\psi$ on
$\mathbb{N}$ such that $\beta(n)=\gamma\psi(n)$ for each $n\in
dom(\beta)$ (of course, the equality assumes that $n\in dom(\psi)$
and $\psi(n)\in dom(\gamma$)).

In the context of algebra and model theory, the transfer of
Computability Theory was initiated in the 1950-s by A.\ Mostowski
\cite{mo52,mo53}, A.V.\ Kuznetsov \cite{ku56,ku58}, and A.\ 
Fr\"ohlich and J.C.\ Shepherdson \cite{fs56}. The subject was
strongly influenced by the work of M.O.\ Rabin \cite{ra60} and
A.I.\ Mal'cev \cite{ma61}. The seminal paper of A.I.\ Mal'cev was
fundamental for the extensive subsequent work in computable
algebra by the Siberian school of algebra and logic. In parallel,
active research in this area was conducted in the West. The
resulting rich theory was summarized in the monographs
\cite{er80,eg99} and the handbook \cite{eg98}.

\begin{defi}\label{cons}
A structure $\mathbb{B}=(B;\sigma)$  of a finite signature $\sigma$
is called constructivizable iff there is a numbering $\beta$ of $B$
such that all signature predicates and functions, and also the
equality predicate,  are $\beta$-computable. Such a numbering
$\beta$ is called a constructivization of $\mathbb{B}$, and the pair
$(\mathbb{B},\beta)$ is called a constructive structure.
\end{defi}

Recall, in particular,  that a binary relation $P$ on $B$ is
$\beta$-computable (resp.\ $\beta$-computably enumerable) if the
corresponding binary relation $P(\beta(m),\beta(n))$ on
$\mathbb{N}$ is computable (resp.\ computably enumerable). In the
case when $\beta$ is a partial numbering, $P$ is called
$\beta$-computably enumerable if there is a computably enumerable
binary relation $\widehat{P}$ on $\mathbb{N}$ such that
$\{(m,n)\mid~P(\beta(m),\beta(n))\}=\widehat{P}\cap(D\times D)$
where $D=dom(\beta)$.

The notion of a constructivizable  structure is equivalent to the
notion of a computably presentable structure popular in the western
literature. Obviously, $(\mathbb{B},\beta)$ is  a constructive
structure iff given a quantifier-free $\sigma$-formula
$\phi(v_1,\ldots,v_k)$ with free variables among $v_1,\ldots,v_k$
and given $n_1,\ldots,n_k\in\mathbb{N}$, one can compute the
truth-value $\phi^\mathbb{B}(\beta(n_1),\ldots,\beta(n_k))$ of
$\phi$ in $\mathbb{B}$ on the elements
$\beta(n_1),\ldots,\beta(n_k)\in B$.

\begin{defi}\label{scons}
A structure $\mathbb{B}=(B;\sigma)$ of a finite signature $\sigma$
is called strongly constructivizable iff there is a numbering
$\beta$ of $B$ such that, given a first-order $\sigma$-formula
$\phi(v_1,\ldots,v_k)$ with free variables among $v_1,\ldots,v_k$
and given $n_1,\ldots,n_k\in\mathbb{N}$, one can compute the
truth-value $\phi^\mathbb{B}(\beta(n_1),\ldots,\beta(n_k))$ of
$\phi$ in $\mathbb{B}$ on the elements
$\beta(n_1),\ldots,\beta(n_k)\in B$. Such a numbering $\beta$ is
called a strong constructivization of $\mathbb{B}$, and the pair
$(\mathbb{B},\beta)$ is called a strongly constructive structure.
\end{defi}

By the definitions above, any strongly  constructivizible structure is
constructivizible and has a decidable first-order theory. Note
that the notion of a strongly constructive structure is equivalent
to the notion of a decidable structure  popular in the
western literature.

We illustrate the introduced notions  by some number structures.
Let
 $\mathbb{N}=(N;<,+,\cdot,0,1)$ be the ordered semiring of
naturals,
 $\mathbb{Z}=(Z;<,+,\cdot,0,1)$  the ordered ring of
integers,
 $\mathbb{Q}=(Q;<,+,\cdot,0,1)$  the ordered field of
rationals,
 $\mathbb{R}=(R;<,+,\cdot,0,1)$  the ordered field
of  reals,
 $\mathbb{R}_c=(R_c;<,+,\cdot,0,1)$  the ordered field
of computable reals \cite{wei},
 and $\mathbb{A}=(A;<,+,\cdot,0,1)$
 the ordered field of algebraic reals  \cite{wa67} (by definition, the algebraic
reals are the real roots of polynomials with rational
coefficients). As is well known \cite{wei}, any algebraic real is
computable, so $\mathbb{A}$ is a substructure of $\mathbb{R}_c$.

As is well-known, the fields $\mathbb{A}$,  $\mathbb{R}_c$ and
$\mathbb{R}$ are real closed (we use some standard algebraic
notions which may be found e.g.\ in \cite{wa67}.)  The following
properties of the mentioned number structures are well-known.
Details and additional  references may be found in the vast
literature on computable rings and fields (see e.g.\ 
\cite{mo66,er68,mn79,er74,st95,st99}).

\begin{exa}\label{exconstr}\leavevmode
  \begin{enumerate}
 \item The  structures $\mathbb{N}$, $\mathbb{Z}$,
$\mathbb{Q}$ are constructivizable but not strongly
constructivizable.

\item The structure $\mathbb{A}$ is strongly  constructivizable.

\item The structure $\mathbb{R}_c$ is not  constructivizable, but
there is a partial numbering $\rho$ of $\mathbb{R}_c$ such that
the field operations are  $\rho$-computable and the relation $<$ is
$\rho$-computably enumerable.
\end{enumerate}
\end{exa}

For this paper, Example~\ref{exconstr}~(2) is of a special interest. We also
need some extensions of this assertion which may be deduced from some known
facts about computable fields (and from classical algebraic facts in
\cite{wa67}).

First we recall definition of the partial numbering $\rho$ of $\mathbb{R}_c$ mentioned in the Example 3. Let
$\varkappa$ be a constructivization of $\mathbb{Q}$ and $\{\varphi_n\}$
be a standard numbering of the computable partial functions on
$\mathbb{N}$. A sequence $\{x_n\}$ in $\mathbb{R}$ is called {\em fast
Cauchy} iff $\forall n\forall i>n(|x_i-x_n|<2^{-n})$. Now, define
$\rho$ as follows: $\rho(n)=x$ iff
$\varphi_n$ is total,  $\{\varkappa\varphi_n(i)\}_i$ is fast
Cauchy, and converges to $x$.
Let us collect some facts relating the introduced notions.

\begin{lem}\label{comreals}\leavevmode
 \begin{enumerate}
 \item Let $\mathbb{B}$ be an ordered subfield of $\mathbb{R}$ and
$\beta$ a constructivization of $\mathbb{B}$.  Then
$\beta\leq\rho$, in particular $\mathbb{B}\subseteq\mathbb{R}_c$.
 \item Let $\mathbb{B}$ be a subfield of $(\mathbb{R};+,\cdot,0,1)$ and
$\beta$ a constructivization of  $\mathbb{B}$ such that
$\beta\leq\rho$. Then $\beta$ is a constructivization of the
ordered field $(\mathbb{B};<)$.
 \item Let $\mathbb{B}$ be a real closed subfield of $(\mathbb{R};+,\cdot,0,1)$ and
$\beta$ a constructivization of $\mathbb{B}$. Then $\beta$ is a
strong constructivization of the ordered field $(\mathbb{B};<)$.
 \end{enumerate}
\end{lem}

\begin{proof}\leavevmode
\begin{enumerate}
  \item Since $\beta$ is a  constructivization,
$\varkappa\leq\beta$. Hence, for some computable functions $f,g$
we have $\varkappa f(n,i)<\beta(n)<\varkappa g(n,i)$ and
$\varkappa g(n,i)-\varkappa f(n,i)<2^{-i}$. Let $h$ be a
computable function satisfying $\varkappa h(n,i)=(\varkappa
g(n,i)-\varkappa f(n,i))/2$. Then $\{\varkappa h(n,i)\}_i$ is a
fast Cauchy sequence converging to $\beta(n)$. Choosing a
computable function $u$ with $h(n,i)=\varphi_{u(n)}(i)$ we see
that $\beta\leq\rho$ via $u$.

\item It suffices to show that $<$ is $\beta$-computable.  Since
$\beta\leq\rho$ and $<$ is $\rho$-computably enumerable, $<$ is
also $\beta$-computably enumerable. Hence, given $n$ one can
compute which of the alternatives
$\beta(n)<0,\beta(n)=0,\beta(n)>0$ holds. Thus, $<$ is
$\beta$-computable.

\item Since $\mathbb{B}$ is real closed,  $0\leq\beta(n)$ is
equivalent to $\exists m(\beta(n)=\beta(m)^2)$. Then $\leq$ and
$<$ are $\beta$-computably enumerable. As in the previous
paragraph, $<$ is $\beta$-computable, hence $\beta$ is a
constructivization of $(\mathbb{B};<)$. By the Tarski quantifier
elimination for real closed fields, given any first order
$\sigma$-formula $\phi$, $\sigma=\{<,+,\cdot,0,1\}$, one can
compute a quantifier-free $\sigma$-formula equivalent to $\phi$ in
$(\mathbb{B};<)$. Thus, $\beta$ is a strong constructivization of
$(\mathbb{B};<)$.\qedhere
\end{enumerate}
\end{proof}

\begin{lem}\label{rfields}
 For any finite set $F\subseteq\mathbb{R}_c$ there is a strongly
constructive real closed ordered subfield $(\mathbb{B},\beta)$ of
$\mathbb{R}_c$ with $F\subseteq B$.
\end{lem}

\proof  If $F\subseteq \mathbb{A}$ we can take
$\mathbb{B}=\mathbb{A}$ and $\beta=\alpha$,  where $\alpha$ is a strong
constructivization of $\mathbb{A}$. Otherwise, let $x$ be the
least element of $F\setminus\mathbb{A}$, so in particular $x$ is a
computable transcendental real number. Let
$\mathbb{D}=\mathbb{A}(x)$ be the subfield of $\mathbb{R}_c$
generated by $\mathbb{A}\cup\{x\}$ and $\delta$ be the numbering
of $\mathbb{D}$ induced by the strong constructivization $\alpha$
of $\mathbb{A}$ and the G\"odel numbering of $\sigma$-terms with
the variable $x$. Since $\alpha\leq\rho$ and $x\in\mathbb{R}_c$,
$\delta\leq\rho$. Moreover, from the well-known structure of
$\mathbb{D}$ it follows that $\delta$ is a constructivization of
$(\mathbb{D};<)$.

Let now $\mathbb{A}_1$ be  the real algebraic closure of
$\mathbb{D}$ in $\mathbb{R}_c$. As is well known (see e.g.\ 
\cite[Theorem~3, p.\ 101]{er74}), $\mathbb{A}_1$ is
constructivizable, even strongly constructivizable by item (3) of
Lemma~\ref{comreals}. If $F\subseteq \mathbb{A}_1$ we can take
$\mathbb{B}=\mathbb{A}_1$ and $\beta=\alpha_1$, where  $\alpha_1$ is a strong
constructivization of $\mathbb{A}_1$. Otherwise, iterate the
construction $\mathbb{A}\mapsto\mathbb{A}_1$ sufficiently many
times in order to get the desired $\mathbb{B}$.
 \qed
 
 \begin{rem}\label{remfields}
 The proof of the last lemma is non-constructive (i.e., from $\rho$-indices of
 elements of $F$ one cannot compute a constructivisation $\beta$). This lemma
 and Theorem~\ref{main} which is based on it are   ``pure existence theorems''.
\end{rem}

Let $(\mathbb{B},\beta)$  be a strongly constructive real closed ordered
subfield of  $\mathbb{R}_c$. Then one can compute, given a
polynomial $p(x)=a_0+a_1x^1\cdots+a_kx^k$ with  coefficients in
$\mathbb{B}$ (i.e., given a string $n_0,\ldots,n_k$ of naturals
with $\beta(n_0)=a_0,\ldots,\beta(n_k)=a_k$) the string
$r_1<\cdots <r_m$, $m\geq 0$, of all distinct real roots of $p(x)$
(i.e., a string $l_1,\ldots,l_m$ of naturals with
$\beta(l_1)=r_1,\ldots,\beta(l_m)=r_m$), as well as the
multiplicity of any root $r_j$. This fact immediately implies

\begin{lem}\label{ceigen}
 Let $(\mathbb{B},\beta)$ be a  strongly constructive real
closed ordered subfield of  $\mathbb{R}_c$. Given a symmetric $n\times
n$-matrix $M$ with coefficients in $\mathbb{B}$, one can compute (w.r.t. $\beta$)
an  orthonormal basis $(\mathbf{v}_1,\ldots,\mathbf{v}_n)\in
\mathbb{B}^n$ of eigenvectors of $M$.
\end{lem}

\proof  Let $(\lambda_1,\ldots,\lambda_n)$ be a string
of all complex roots of the characteristic polynomial ${\rm det}(\lambda
I_n-M)$ taken with their multiplicities. Since $M$ is symmetric,
$\lambda_1,\ldots,\lambda_n\in \mathbb{R}$. Since $\mathbb{B}$ is
real closed and the coefficients of  the characteristic polynomial
are  in $\mathbb{B}$, $\lambda_1,\ldots,\lambda_n\in \mathbb{B}$.
Since $(\mathbb{B},\beta)$ is strongly constructive, given
$\beta$-names for the coefficients of $M$ one can compute
$\beta$-names for $\lambda_1,\ldots,\lambda_n$ (without loss of
generality we may even assume that
$\lambda_1\leq\cdots\leq\lambda_n$). Repeating well-known
computations from linear algebra (cf.\ e.g.\ the proof of \cite[Theorem~13]{zb04}) one can compute the desired eigenvectors
$\mathbf{v}_1,\ldots,\mathbf{v}_n\in \mathbb{B}^n$ such that
$M\cdot\mathbf{v}_i=\lambda_i\cdot\mathbf{v}_i$ for each
$i=1,\ldots,n$.
 \qed

\begin{rem}\label{ceigen1}
 Of course, the orthonormal basis of eigenvectors is not unique. It is only important that some such basis is computable (w.r.t. $\beta$).
\end{rem}

\subsection{Computability on metric spaces}\label{compu}

We use the TTE-approach to computability over metric spaces
developed in the K. Weihrauch's school (for more details see e.g.\ 
\cite{wei,wz02,br03,bhw} and references therein). Let $(M,d)$ be a
metric space. A  sequence $\{x_n\}$ in $M$ is called {\em fast
Cauchy} iff $d(x_i,x_n)<2^{-n}$ for all $n$ and $i>n$. The
following lemma is straightforward.

\begin{lem}\label{fcauchy}
Let $(M,d)$ be a metric space and let $x,x_n,x_{n,m}\in M$ for all $m,n\in\mathbb{N}$.
\begin{enumerate}
 \item If $\{x_n\}$ is fast Cauchy and converges to $x$ then
 $\forall n(d(x,x_n)\leq 2^{-n})$.
 \item If $\forall n(d(x,x_n)\leq 2^{-n})$ then $\{x_n\}$ converges to $x$
 and $\{x_{n+1}\}$ is fast Cauchy.
 \item Let for any $n$ $\{x_{n,m}\}_m$ is fast Cauchy and converges to $x_n$,
 and let  $\{x_n\}$ is fast Cauchy and converges to $x$. Then $\{x_{n+2,n+2}\}$
 is fast Cauchy and converges to $x$.
\end{enumerate}
\end{lem}

Let $ \mathcal{N}=\omega^\omega$ be the Baire  space (instead of
the Baire space people often use in this context the Cantor space $\Sigma^\omega$
of infinite words over a finite alphabet $\Sigma$ containing at
least two symbols.). Relate to any function $\nu:\mathbb{N}\to M$
the partial function $\nu^\ast$ from $\mathcal{N}$ to $M$ as
follows: $\nu^\ast(p)=x$ iff the sequence $\{\nu_{p(n)}\}$ is
fast Cauchy and converges to $x$.

\begin{lem}\label{fcauchy1}
Let $(M,d)$ be a metric space and  $\mu,\nu:\mathbb{N}\to M$ be
such that $\mu\leq\nu^\ast$ (i.e., $\mu=\nu^\ast\circ f$ for a
computable function $f:\mathbb{N}\to \mathcal{N}$). Then
$\mu^\ast\leq\nu^\ast$ (i.e., $\mu^\ast=\nu^\ast\circ g$ for a
computable partial function $g$ on $\mathcal{N}$).
\end{lem}

\proof  For any $p\in dom(\mu^\ast)$,
$\mu^\ast(p)=lim_n\mu_{p(n)}$  and $\{\mu_{p(n)}\}$ is fast
Cauchy. For each $n$, $\mu(p(n))=\rm{lim}_m\nu_{f(p(n))(m)}$ and
$\{\nu_{f(p(n))(m)}\}_m$ is fast Cauchy. By item (3) of Lemma~\ref{fcauchy},
$\{\nu_{f(p(n+2))(n+2)}\}_n$ is fast Cauchy and
converges to $\mu^\ast(p)$. Let $g$ be the computable function on
$\mathcal{N}$ defined by $g(p)(n)=f(p(n+2))(n+2)$. Then $g$ has
the desired property.
 \qed

\begin{defi}\label{cmetric}
 A computable metric space is a triple
$(M,d,\mu)$ where $(M,d)$ is a metric space and
$\mu:\omega\rightarrow M$ is a numbering of a dense subset ${\rm
rng}(\mu)$ of $M$ such that $\{d(\mu_m,\mu_n)\}$ is a computable
double  sequence of reals. The partial surjection $\mu^\ast$ from
$\mathcal{N}$ onto $M$ is called the Cauchy representation of
$(M,d,\mu)$.
\end{defi}

Note that the computability of the double sequence
$\{d(\mu_m,\mu_n)\}$ is equivalent to the computable enumerability
of the set
 $$
\{(i,j,q,r)\mid i,j\in\omega,\,q,r\in{\mathbb
Q},\,q<d(\nu_i,\nu_j)<r\}.
 $$

A partial function $f:M\rightharpoonup M_1$ on the elements of
computable metric spaces $(M,d,\nu)$ and $(M_1,d_1,\nu_1)$ is {\em
computable} if there is a computable partial function $\hat{f}$
on $\mathcal{N}$ which
realizes $f$ w.r.t.  the Cauchy representations of $M$ and
$M_1$, i.e., $\nu^\ast_1(\hat{f}(p))=f(\nu^\ast(p))$ for each
$p\in \operatorname{dom}(\nu^\ast)$ (in other words, if $\{\nu (p(n))\}$ is a fast Cauchy sequence converging to $x\in M$ then $\{\nu_1(\hat{f}(p)(n))\}$ is a fast Cauchy sequence converging to $f(x)\in M_1$).

A standard  example of a computable metric space is
$(\mathbb{R},d,\varkappa)$ where $d(x,y)=|x-y|$ is the standard
metric on $\mathbb{R}$ and $\varkappa$ is a constructivization of
$\mathbb{Q}$ (see the previous subsection). A less standard
example is $(\mathbb{R},d,\beta)$ where $\beta$ is a strong
constructivization  of a real closed ordered subfield $\mathbb{B}$
of $\mathbb{R}_c$.  Though formally different, these two
computable metric spaces are equivalent in the following sense:

\begin{lem}\label{cmetric1}
 The Cauchy representations $\varkappa^\ast,\beta^\ast$ of
$\mathbb{R}$ induced by the numberings $\varkappa,\beta$
respectively, are equivalent, i.e.\ $\varkappa^\ast\leq\beta^\ast$
and $\beta^\ast\leq\varkappa^\ast$.
\end{lem}

\proof   Since
$\varkappa\leq\beta$ and $\beta\leq\beta^\ast$, we have
$\varkappa\leq\beta^\ast$, hence $\varkappa^\ast\leq\beta^\ast$
by Lemma~\ref{fcauchy1}.  For the converse reduction, by Lemma~\ref{fcauchy1}
it suffices to show that $\beta\leq\varkappa^\ast$.
This follows from item~(1) of Lemma~\ref{comreals}.
 \qed

\subsection{Spaces under consideration}\label{spaces}

For any $n\geq 1$, the vector space  ${\mathbb R}^n$ carries the
sup-norm $||x||_\infty=\operatorname{max}\{|x_i|\}$  and the Euclidean norm
$||x||_2=\sqrt{\sum x_i^2}$; we denote the corresponding metrics
by $d_\infty$ and $d_2$, respectively.

Define the function $\varkappa^n:\mathbb{N}\to{\mathbb R}^n$ by
$\varkappa^n\langle
k_1,\ldots,k_n\rangle=(\varkappa(k_1),\cdots,\varkappa(k_n))$
where $\langle\cdot\rangle$ is a computable coding function of
$n$-tuples of naturals and $\varkappa$ is a constructivization of
$\mathbb{Q}$ (see Section~\ref{compcount}).  Let
$(\mathbb{B},\beta)$ be a strongly constructive real closed ordered
subfield of  $\mathbb{R}_c$.  Define the function
$\beta^n:\mathbb{N}\to{\mathbb R}^n$ in the same way
as $\varkappa^n$, with $\varkappa$ replaced by  $\beta$.

\begin{lem}\label{rnspace}
For any $n\geq 1$ and $d\in\{d_\infty,d_2\}$,  $({\mathbb
R}^n,d,\varkappa^n)$ and $({\mathbb R}^n,d,\beta^n)$ are
equivalent computable metric spaces.
\end{lem}

\proof  For $n=1$ this follows from Lemma~\ref{cmetric1}
because $d_\infty=d_2$. For $n\geq 2$,   computability of the
spaces and the reducibility $\varkappa^n\leq\beta^n$ are obvious.
Since $\beta\leq\rho$ by item~(1) of Lemma~\ref{comreals}, there is
a computable function $f$ on ${\mathbb N}$ such that
$d_\infty(\beta(k),\varkappa f\langle k,l\rangle)\leq 2^{-l}$ for all
$k,l$. By Lemma~\ref{fcauchy1} and the argument of Lemma~\ref{cmetric1},
for the metric $d=d_\infty$ it suffices to find a
computable function $g$ on ${\mathbb N}$ such that
$d_\infty(\beta^n(k),\varkappa^ng\langle k,l\rangle)\leq 2^{-l}$ for
all $k,l$. Define $g$ by $g\langle k,l\rangle=\langle f\langle
k_1,l\rangle,\ldots,f\langle k_n,l\rangle\rangle$ where $k=\langle
k_1,\ldots,k_n\rangle$. Then we have
 $$d_\infty(\beta^n(k),\varkappa^ng\langle k,l\rangle)=
 \operatorname{sup}\{ d_\infty(\beta(k_1),\varkappa f\langle k_1,l\rangle),\ldots,
 d_\infty(\beta(k_n),\varkappa f\langle k,l\rangle)\}\leq 2^{-l}
 $$
 which completes the case $d=d_\infty$.

For $d=d_2$, the assertion follows from the   obvious 
estimate  $d_2(x,y)\leq \sqrt{n} d_\infty(x,y)$. 
 \qed

We will consider some subspaces of the introduced metric spaces,
in particular the space $S\subseteq{\mathbb R}^{n\times n}$ of
symmetric real matrices, the space $S_{+}$  of
symmetric real positively definite matrices, and the
$m$-dimensional unitary cube $Q=[0,1]^m$. For these subspaces the
analog of Lemma~\ref{rnspace} clearly holds.

In the study of difference equations,  some norms on the spaces of
grid functions are quite useful. For any $N\geq 0$, let $G=G_N$ be
the {\it uniform grid on} $Q$ formed by the binary-rational
vectors $(x_1,\ldots,x_m)$ where $x_i=\frac{y_i} {2^N}$ and
$y_i\in\{0,1,\ldots,2^N\}$. Note that the number of such vectors
is $(2^{N}+1)^m$, so the set $\mathbb{R}^{G}$ of grid functions
$f:G_N\to\mathbb{R}$ may be identified with $\mathbb{R}^{
(2^{N}+1)^m}$ while the set $(\mathbb R^n)^{G}$ of grid functions
$f:G_N\to\mathbb{R}^n$ may be identified with $\mathbb{R}^{n\cdot
(2^{N}+1)^m}$. In the last case, we obtain the following norms
 $$||\varphi||_s=\operatorname{max}_{x\in G_N}||\varphi(x)||, \;
 ||\varphi||^2_{L_2}=h^m\sum_{x\in G_N}\langle\varphi(x),\varphi(x)\rangle.$$
 Note that $d_s$ coincides with $d_\infty$ for the corresponding metric spaces, and Lemma~\ref{rnspace} applies to these spaces.

For all rational $\tau>0$ and integers $N\geq0$ and $L\geq1$, let
$G_N^\tau$ be the grid in $Q\times[0,T]$, $T=L\tau$, with step
$h=\frac{1} {2^N}$ on the space coordinates $x_i$ and step $\tau$
on the time coordinate $t$. Just as above, we can define the sup- and $L_2$-norms on the vector space  $M=\mathbb{R}^{G^\tau_N}$ of grid
functions on such grids   and Lemma~\ref{rnspace} applies to the corresponding metric spaces.

 The vector space $M$ carries also another natural norm (called the
$sL_2$-norm.) The $sL_2$-norm of a grid function
$f:G_N^\tau\rightarrow{\mathbb R}$ is defined by
 $$||f||^2_{sL_2}=\max_{0\leq l\tau\leq T}\left(h^m\sum_{x\in
G_N} f^2(x,l\tau)\right).$$
 Let $d_{sL_2}$ be the corresponding
metric on $M$. Let $(\mathbb{B},\beta)$ be  a strongly
constructive real closed ordered subfield of  $\mathbb{R}_c$.  Let $\mu$
(resp.\ $\nu$) be the numbering of the set $\{f\mid
f:G_N^\tau\rightarrow{\mathbb Q}\}$ (resp.\ $\{f\mid
f:G_N^\tau\rightarrow{\mathbb B}\}$) induced by the natural
numbering $G_0^\tau,G_1^\tau,\cdots$ of all such grids and the
constructivization $\varkappa$ of ${\mathbb Q}$ (resp.\ the strong
constructivization $\beta$ of ${\mathbb B}$). The following analog
of Lemma~\ref{rnspace} is an easy corollary of the estimate
$d_{sL_2}(f,g)\leq d_s(f,g)$ which follows from the definition of
$sL_2$-norm.

\begin{lem}\label{sl_2space}
In the notation of the previous paragraph, $(M,d_{sL_2},\mu)$  and $(M,d_{sL_2},\nu)$ are equivalent computable metric spaces. This extends in the obvious way to the space $(\mathbb{R}^n)^{G^\tau_N}$ of  grid functions
$f:G^\tau_N\to\mathbb{R}^n$.
\end{lem}

We will work with several functional spaces most of which are
subsets of the set $C({\mathbb R}^m,{\mathbb R}^n)\simeq
C({\mathbb R}^m,{\mathbb R})^n$ of integrable continuous functions
$\varphi:{\mathbb R}^m\rightarrow{\mathbb R}^n$ equipped with  the
$L_2$-norm. In particular, we deal with the space $C(Q,{\mathbb
R}^n)\simeq C(Q,{\mathbb R})^n$ (resp.\ $C^k(Q,{\mathbb R}^n)$) of
continuous (resp.\ $k$-time  continuously differentiable) functions
$\varphi:Q\rightarrow{\mathbb R}^n$ equipped with the $L_2$-norm
 $$||\varphi||_{L_2}=\left(\int_Q|\varphi(x)|^2dx)\right)^{\frac{1}{2}},\; |\varphi(x)|^2=
\langle\varphi,\varphi\rangle=\sum\limits_{i=1}^n\varphi^2_i(x).$$
 We will also use the  sup-norm $$||\varphi||_s=\sup_{x\in
Q}|\varphi(x)|,\ ||f||_s=\sup_{(x,t)\in Q\times[0,T]}|f(x,t)|$$ on
$C(Q,\mathbb R^n)$ or $C(Q\times[0,T],\mathbb R^n)$ and the
$sL_2$-norm
 $$||u||_{sL_2}=\sup_{0\leq t_0\leq
 T}\sqrt{\int_Q|u(x,t_0)|^2dx}$$
  on
$C(Q\times[0,T],{\mathbb R^n})$ where $T>0$.  Whenever we want to
emphasize the norm we use notations like $C_{L_2}(Q,{\mathbb
R}^n)$, $C_s(Q,{\mathbb R}^n)$ or $C_{sL_2}(Q\times[0,T],{\mathbb
R}^n)$.

Associate to any grid function $f_N:G_N\rightarrow{\mathbb Q}$ the
continuous extension $\tilde{f}_N:Q\rightarrow{\mathbb R}$ of $f$
obtained by the piecewise-linear interpolation on each coordinate. Such
interpolations known also as {\it multilinear interpolations} are
the simplest class of splines  see e.g.\ \cite{sz59,so74,za,ba86}).
Note that the restriction of $\tilde{f}_N$ to any grid cell is a
polynomial of degree $m$, see an example in Subsection~\ref{spaces}. The
extensions $\tilde{f}_N$ induce a countable dense set in
$C(Q,{\mathbb R}^n)$ (or $C(Q\times[0,T],{\mathbb R^n})$) with any
of the three norms.

Let again $(\mathbb{B},\beta)$  be a strongly constructive  real
closed ordered subfield of  $\mathbb{R}_c$. Define
$\tilde{\beta},\tilde{\varkappa}:\mathbb{N}\to C(Q,{\mathbb R}^n)$
by $\tilde{\beta}\langle N,l\rangle=\widetilde{\beta_N^p(l)}$ where
$p$ is the number of grid points in $G_N$ and $\beta_N^p$ is the
numbering of grid functions $f:G_N\to\mathbb{B}^n$
($\tilde{\varkappa}$ is defined similarly). Define also
$\tilde{\mu},\tilde{\nu}:\mathbb{N}\to C(Q\times[0,T],{\mathbb
R}^n)$ in the same way, starting from the numberings $\mu,\nu$ above
and the natural numbering of all the grids $G^\tau_N$ with rational
positive $\tau$. The next fact follows from Lemmas~\ref{rnspace},
\ref{sl_2space} and the well-known estimates
$||\tilde{f}||\leq||f||$ where $||\cdot||$ is any of the three norms
(see e.g.\ \cite[p.\ 187--189]{ba86}, \cite[p.\ 335]{sz59}).

\begin{lem}\label{fspace}\leavevmode
 \begin{enumerate}
 \item For any $n\geq 1$ and $d\in\{d_s,d_{L_2}\}$,
$(C(Q,{\mathbb R}^n),d,\tilde{\varkappa})$ and $(C(Q,{\mathbb
R}^n),d,\tilde{\beta})$ are equivalent computable metric spaces.
 \item In the notation before the formulation of lemma,
$(C(Q\times[0,T],{\mathbb R}^n),d_{sL_2},\tilde{\mu})$ and
$(C(Q\times[0,T],{\mathbb R}^n),d_{sL_2},\tilde{\nu})$ are equivalent
computable metric spaces.
 \end{enumerate}
\end{lem}

Let again $G$ be the grid in $Q$ with step $h=\frac{1}{2^N}$ on each
coordinate. From well-known facts of Computable Analysis
\cite{wei} it follows that the restriction
$\varphi\mapsto\varphi|_G$ is a computable operator from
$C_s(Q,{\mathbb R}^n)$ to $({\mathbb R}^n)^G$. From well-known
properties of the multilinear interpolations (see e.g.\ 
\cite{go71,za}) it follows that $f\mapsto\tilde{f}$ is a
computable operator from $(({\mathbb R}^n)^G)_s$ to
$C_{L_2}(Q,{\mathbb R}^n)$ (see also the estimate~\eqref{14}
below).

Along with the mentioned norms,  we sometimes use their
$A$-modifications, for a given matrix $A$. In particular, the
$A$-modification of the $L_2$-norms is defined by
 $$
||\varphi||_{A,L_2}=\sqrt{\int_Q\langle A\varphi,\varphi\rangle dx}
 $$
 while the $A$-modification of the $sL_2$-norms is defined by
 $$
||u||_{A,sL_2}=\sup_{0\leq t_0\leq
 T}\sqrt{\int_Q\langle Au(x,t_0),u(x,t_0)\rangle dx},
 $$
 and in a similar way for the grid norms.

\section{Finite-dimensional approximations}\label{fdim}

In this section we describe the construction of difference
operators approximating  the differential problem considered in
this paper and establish their basic properties. Subsection~\ref{basicdif}
recalls some relevant notions and general facts
about difference schemes. In Subsection~\ref{condif} we describe
the difference scheme \cite{go76} for the symmetric hyperbolic systems under
consideration. In Subsection~\ref{difop} we establish some basic
properties of the corresponding difference operators.

\subsection{Basic facts about difference schemes}\label{basicdif}

Here we briefly recall some relevant notions and facts about
difference schemes (for more details see any book on the subject,
e.g.\ \cite{gr,jo66,strik,tre96}).

Let us consider the boundary-value problem~\eqref{sist} for a
(system of) PDEs. Difference approximations to~\eqref{sist}  are
written in the form
  \begin{equation}\label{pdeappr}
L_h{\bf u}^{(h)}={\bf f}^{(h)},\;{\mathcal L}_h{\bf
u}^{(h)}=\varphi^{(h)}
 \end{equation}
where $L_h,{\mathcal L}_h$ are difference operators (which are in
our case linear), and all functions are defined on some grids in
$\Omega$ or $\Gamma\subseteq\partial\Omega$ (the grids are not
always uniform, as in our simplest case). For simplicity we will use
the restriction notation ${\bf g}|_{G_k}$ to denote the restriction
of ${\bf g}:\Omega\rightarrow\mathbb{R}^n$ to the grid $G_k$ in
$\Omega$ though in general the restriction operator may be more
complicated. Both sides of~\eqref{pdeappr} depend on the grid step
$h$.

Note that  in this paper we consider a little more complicated
grids than the uniform grids discussed above, namely grids with
the integer time steps $l\tau$, $l\geq 0$, (for some $\tau>0$) and
half-integer steps $x_{i+\frac{1}{2}}=(i+\frac{1}{2})h$ for the space variables.
The theory for such slightly modified grids remains  the same.

Let the space of grid functions defined on the same grid as ${\bf
f}^{(h)}$ (resp.\ as ${\bf u}^{(h)}$, $\varphi^{(h)}$) carry some
norm $||\cdot||_{F_h}$ (resp.\ some norms $||\cdot||_{U_h}$,
$||\cdot||_{\Phi_h}$). Note that in our case these will be
$L_2$ and $sL_2$-norms defined in Section~\ref{prel}.

\begin{defi}\label{appro}
Difference equations~\eqref{pdeappr}, also called difference
schemes, {\em approximate the differential equation~\eqref{sist}
with order of accuracy $l$ (where $l$ is a positive integer) on a
solution ${\bf u}({\bf x},t)$ of~\eqref{sist}} if
\begin{eqnarray*}
 ||(L{\bf
u})|_{G_k}-L_h{\bf u}^{(h)}||_{F_h}\leq M_1h^l,\;
||f|_{G_k}-f^{(h)}||_{F_h}\leq M_2h^l,\\ ||({\mathcal L}{\bf
u})|_{G_k}-{\mathcal L}_h{\bf u}^{(h)}||_{\Phi_h}\leq M_3h^l \text{
and
 }||\varphi|_{G_k}-\varphi^{(h)}||_{\Phi_h}\leq M_4h^l
\end{eqnarray*}
 for some
constants $M_1,M_2,M_3$ and $M_4$ not depending on $h$ and
$\tau$.
\end{defi}

The definition is usually checked by working with the Taylor
series for the corresponding functions, thus the constants $M_i$
depend on the derivatives of the functions $u$ and $f$. As a result,
the degrees of smoothness of the functions become essential when
one is interested in the order of accuracy of a difference scheme.
Note that the definition assumes the existence of a solution of~\eqref{sist}.
For the problems~\eqref{sist_1} and~\eqref{sist_2}
it is well-known (see e.g.\ \cite{fri,go71,mi,evans}) that there is a
unique solution.

The following notion identifies a property of difference schemes
which is crucial for computing ``good'' approximations to the
solutions of~\eqref{sist}.

\begin{defi}\label{stable}
Difference scheme~\eqref{pdeappr}  is called {\em stable} if its
solution ${\bf u}^{(h)}$ satisfies
 $$||{\bf u}^{(h)}||_{U_h}\leq
N_1||f^{(h)}||_{F_h}+N_2||\varphi^{(h)}||_{\Phi_h}$$
 for some
constants $N_1$ and $N_2$ not depending on $h$, $\tau$, $f^{(h)}$
and $\varphi^{(h)}$.
\end{defi}

For non-stationary processes (depending explicitly on the time
variable $t$, as~\eqref{sist_1},~\eqref{sist_2}), the difference
equation~\eqref{pdeappr} may be rewritten in the equivalent
recurrent form ${\bf u}^{[l+1]}=R_h{\bf u}^{[l]}+\tau\rho^{[l]}$
where ${\bf u}^{[0]}$ is known, ${\bf u}^{[l]}$ is the restriction
of the solution to the time level $t=l\tau$, $l\geq0$,
$\rho^{[l]}$ depends only on $f$ and $\varphi$, $R_h$ is the
difference operator obtained from $L_h$ in a natural way. It is
known (see e.g.\ \cite{gr}) that the stability of~\eqref{pdeappr}
on the interval $0<t<T$ is equivalent to the uniform boundedness
of the operators $R_h$ and their powers: $||R^m_h||<K$,
$m=1,2,\ldots,\frac{T}{\tau}$, for some constant $K$ not depending
on $h$ and $\tau$. In general, the investigation of the stability
of difference schemes is  a hard task. The most popular tool is
the so called Fourier method \cite{gr,go76,tre96}; for  problems (2), (3) and for the
scheme from the next subsection this is done by using the discrete
energy integral technique in \cite[p.\ 79]{go76}. We will briefly
describe this idea below.

Our main results on the  computability  of solution operators for~
\eqref{sist_1} and~\eqref{sist_2} make an essential use of the
following basic fact from the theory of difference schemes (see
e.g.\ \cite[p.\ 172]{gr},):

\begin{thm}\label{convds}
Let the difference scheme~\eqref{pdeappr} be stable and
approximate~\eqref{sist} on the solution ${\bf u}$ with order $l$.
Then the solution of~\eqref{pdeappr} uniformly converges to the
solution ${\bf u}$  in the sense that $||{\bf u}|_{G^\tau_k}-{\bf
u}^{(h)}||_{U_h}\leq Nh^l$ for some constant $N$ not depending on
$h$ and $\tau$.
\end{thm}

\subsection{Constructing a difference
scheme for symmetric hyperbolic systems}\label{condif}

The difference scheme for the boundary-value problem~
\eqref{sist_2} and the Cauchy problem~\eqref{sist_1} may be
chosen in various ways. The scheme we use is taken from \cite{go76}. It
can be applied to a broader class of systems, including  some systems of
nonlinear equations. We describe it in few stages, letting for
simplicity the righthand part to be zero: $f=0$.

\begin{enumerate}[label=\arabic*.]
\item First we describe some discretization details. To simplify
notation, we stick to the 2-dimensional case $x_1=x$, $x_2=y$,
$B_1=B$, $B_2=C$, i.e., $m=2$. For $m\geq3$ the
difference scheme is obtained in the same way as for $m=2$ but the
step from $m=1$ to $m=2$ is nontrivial.

Consider the uniform rectangular  grid $G$ on $Q=[0,1]^2$ defined
by the family of lines $\{x=x_i\},\{y=y_j\}$ where $0\leq
i,j\leq2^N$ for some natural number $N$. Let
$h=x_i-x_{i-1}=y_j-y_{j-1}=1/2^N$ be the step of the grid.
Associate to any function $g\in\{u_1,\ldots,u_n\}$ and any fixed
time point $t=l\tau,$ $l\in\mathbb N$, the vector of dimension
$2^{2N}$ with the components
 \begin{equation}\label{init}
g_{i-\frac{1}{2},j-\frac{1}{2}}=g
\left(\frac{i-\frac{1}{2}}{2^N},\frac{j-\frac{1}{2}}{2^N},t\right)
 \end{equation}
 equal to the values of $g$
in the centers of grid cells, and denote, as in the previous
subsection,
 $${\bf u}^{(h)}=\{{\bf
u}_{i-\frac{1}{2},j-\frac{1}{2}}\}_{1\leq
i,j\leq 2^N, t=l\tau,l=1,\ldots,M}.
 $$
 The initial grid
function, from which the iteration process starts, will be denoted
as ${\bf \varphi}^{(h)}$, which equals to ${\bf u}^{(h)}$
restricted to the time level $t=0$.

Note that, strictly speaking, we work with modifications of the
grids $G_k$ in Subsection~\ref{spaces} when the centers of grid cells are
taken as nodes of the modified grids.

\item Consider the following two auxiliary one-dimensional systems
with parameters obtained by fixing any of the variables $x,y$:
 \begin{equation}\label{auxil}
 A\frac{\partial {\bf u}}{\partial t}+B\frac{\partial
{\bf u}} {\partial x}=0,\;A\frac{\partial {\bf u}}{\partial
t}+C\frac{\partial {\bf u}} {\partial y}=0,
 \end{equation}
 where ${\bf u}=(u_1,u_2,\ldots,u_n)^T$. Transform the systems
 into their canonical forms
\begin{equation}\label{canon}
 \frac{\partial {\bf v}_x}{\partial
t}+M_x\frac{\partial {\bf v}_x} {\partial x}=0,\ \
\frac{\partial {\bf v}_y}{\partial t}+M_y\frac{\partial {\bf v}_y}
{\partial y}=0
 \end{equation}
 via the linear transformations ${\bf v}_x=T_x^{-1}{\bf u}$ and
${\bf v}_y=T_y^{-1}{\bf u}$ defined as follows (see
\cite{go76,ss09} for additional details):
 \begin{equation}\label{T}T_x=LDK_x,\ {\bf u}=T_x{\bf v}_x,\end{equation}
 and in a similar way on the $y$-coordinate. Here the orthogonal matrix $L$ transforms  the
matrix $A$ to its canonical form
 $L^*AL=\Lambda={\rm
 diag}\{\lambda_1,\lambda_2,\ldots,\lambda_n\}$,
$D=\Lambda^{-\frac{1}{2}}$. 
The orthogonal matrix 
  $K_x$ transforms  the matrix $D^*L^*BLD$ to its diagonal form $M_x$. And similarly for the second auxiliary system in~\eqref{auxil}.

Note that the matrices $L,K_x,K_y$ consisting  of eigenvectors of the
corresponding symmetric matrices are not uniquely  defined. We
choose some orthonormal eigenvectors and keep them fixed for all
iteration steps. The components of the vectors ${\bf v}_x$, ${\bf
v}_y$ in~\eqref{canon} are called {\em Riemannian invariants};
they are invariant along the characteristics of the corresponding one-dimensional systems, see e.g.\ \cite{go71,go76,evans,mi}.  The existence of
these invariants is a corollary of the hyperbolicity property.

\item Any of the systems~\eqref{canon} in the canonical form consists
of $n$ independent equations of the form
 \begin{equation}\label{eq_1} \frac{\partial w}{\partial
t}+\mu\frac{\partial w} {\partial x}=0,
 \end{equation}
 where $w=w(x,t)$ is a scalar function and $\mu\in\mathbb R$.
Consider for the  equation~\eqref{eq_1} the following difference
scheme. The function $w(t_0,x)$,  already computed at time level
$t=t_0$ (initially $t=0$; the values on this level are taken from the initial conditions), is substituted by the
piecewise-constant function with the values $w_{i-\frac{1}{2}}$
within the corresponding grid cell $x_{i-1}<x\leq x_i$. Define for each $1\leq i\leq 2^N-1$ auxiliary ``interior'' 
values (called ``large values'' in \cite{go76})  as follows:
 \begin{equation}\label{flux_1}{\mathcal W}_i=
 \begin{cases}w_{i-\frac{1}{2}},\textnormal{
if }\mu\geq0,\\ w_{i+\frac{1}{2}},\textnormal{ if }\mu<0.
 \end{cases}
 \end{equation}
 
In the case of the Cauchy
problem~\eqref{sist_1}, for the auxiliary ``boundary'' values ${\mathcal W}_0$ and ${\mathcal W}_{2^N}$  we use the same formula~\eqref{flux_1} where 
$w_{-\frac{1}{2}}=w_{2^N-\frac{1}{2}}$ and
$w_{2^N+\frac{1}{2}}=w_{\frac{1}{2}}$. 
 The case of the boundary-value problem is described in the step 4 below. 
 
  The values $\{w^{i-\frac{1}{2}}\}$ on the next time level $t=t_0+\tau$ ($\tau$ is a time step
depending on $h$ as specified in the next subsection) are then computed as
 \begin{equation}\label{scheme_1}
w^{i-\frac{1}{2}}=w_{i-\frac{1}{2}}-\mu\frac{\tau}{h}({\mathcal
W}_i- {\mathcal W}_{i-1}).
  \end{equation}
  

Taking the scheme~\eqref{scheme_1} for each equation of the
systems~\eqref{canon}, we obtain for them schemes of the following vector
form:
 $$
  \frac{{\bf
v_x}^{i-\frac{1}{2}}-{\bf
v_x}_{i-\frac{1}{2}}}{\tau}+M_x\frac{({\mathcal V}_x)_i-({\mathcal
V}_x)_{i-1}}{h}=0.
 $$

\item For the boundary-value problem~\eqref{sist_2}, we
compute the vector boundary values $(\mathcal{V}_x)_0,(\mathcal{V}_x)_{2^N}$
with the help of the boundary conditions. On the left boundary
$x=0$  we calculate $m_{+}$ components of $(\mathcal{V}_x)_0$,
corresponding to the positive eigenvalues of the matrix $A^{-1}B$,
from the system of linear equations $\Phi_1^{(1)}(T_x\mathcal{V}_x)_0=0$;
for $m_{-}$ components of $(\mathcal{V}_x)_0$, corresponding to the
negative eigenvalues, we let $(\mathcal{V}_x)_0:=({\bf
v_x})_{\frac{1}{2}}$. The
components corresponding to the zero eigenvalues of $A^{-1}B$ can
be chosen arbitrary since they are multiplied by zero in the
scheme. The values on the right boundary and on both boundaries by
the $y$-coordinate are calculated in a similar way.

\item  Finally, for finding the values $\{{\bf
u}^{i-\frac{1}{2},j-\frac{1}{2}}\}$ on the next time step, we use the system of linear equations
 \begin{equation}\label{scheme_2_sist} A\frac{{\bf
u}^{i-\frac{1}{2},j-\frac{1}{2}}-{\bf
u}_{i-\frac{1}{2},j-\frac{1}{2}}}{\tau}+B\frac{{\mathcal
U}_{i,j-\frac{1}{2}}-{\mathcal
U}_{i-1,j-\frac{1}{2}}}{h}+C\frac{{\mathcal
U}_{i-\frac{1}{2},j}-{\mathcal U}_{i-\frac{1}{2},j-1}}{h}=0
\end{equation}
where
 ${\mathcal U}_{i,j-\frac{1}{2}}=T_x({\mathcal V}_x)_i$ and ${\mathcal
U}_{i-\frac{1}{2},j}=T_y({\mathcal V}_y)_j$ are obtained by applying the  transformations inverse to~\eqref{T}.
\end{enumerate}

\noindent The scheme~\eqref{scheme_2_sist}  was invented by S.K. Godunov; it is described and analysed in all details in \cite{go76}, see also e.g.\ \cite{kps,gv}. It approximates the system~\eqref{sist_1} or~
\eqref{sist_2} with the first order of accuracy  (the
proof of the approximation property is done by means of the Taylor
decomposition).

\subsection{Properties of the difference operators}\label{difop}

Here we establish some  properties of the difference operators from
the previous subsection. 

\begin{lem}\label{complin}\leavevmode
 \begin{enumerate}
 \item The difference operators $\left\{{\bf u}_{i-\frac{1}{2},
j-\frac{1}{2}}\right\}\mapsto\left\{{\bf
u}^{i-\frac{1}{2},j-\frac{1}{2}}\right\}$ and
$\varphi^{(h)}\mapsto u^{(h)}$ are linear.
 \item Let $(\mathbb{B},\beta)$ be  a strongly constructive
real closed ordered subfield of  $\mathbb{R}_c$. Given symmetric matrices
$A, B_1,\ldots,B_m$ with coefficients in $\mathbb{B}$, an initial
grid function $\varphi^{(h)}$ on $G_k$ with  values in
$\mathbb{B}$,  and a number $\tau\in\mathbb{B}$, one can compute
(w.r.t.   $\beta$) the grid function
$ u^{(h)}$ on $G^\tau_k$ (which again has its values in
$\mathbb{B}$).
 \item For any positive $\tau\leq\left(\frac{1}{\tau_x}+\frac{1}{\tau_y}\right)^{-1}$,
where \begin{equation}\label{stab}\tau_x=\max_i\{|\mu_i|:\text{det}(\mu_i A-B)=0\}\cdot h\ \text{ and }\
\tau_y=\max_i\{|\mu_i|: \text{det}(\mu_i A-C)=0\}\cdot h\end{equation}  ($i=1,2,\ldots,n$), the
difference scheme~\eqref{scheme_2_sist} is stable in the sense of
Definition~\ref{stable}.
  \end{enumerate}
\end{lem}

\begin{proof}\leavevmode
\begin{enumerate}
\item The assertion follows from the fact that~
\eqref{scheme_2_sist} is a linear system of equations and the
operator computing the ``large'' values $\mathcal{U}_i$ is linear provided
that the eigenvectors found at step 2 are fixed.

\item The operator $(A,B_1,\ldots,B_m,\varphi^{(h)},\tau)\mapsto
u^{(h)}$  involves only algebraic (iterative) computations,
including the solution of  linear systems of equations with
(previously computed) coefficients in $\mathbb{B}$,  finding of
eigenvalues and eigenvectors of symmetric matrices with
(previously computed) coefficients in $\mathbb{B}$, and comparing
(previously computed) numbers in $\mathbb{B}$, in particular in
the branching operators~\eqref{flux_1},~\eqref{scheme_1}.
Therefore, the assertion follows from Lemma~\ref{ceigen} and other
remarks in Section~\ref{compcount}.

\item This assertion is a standard fact, we give some details for the reader not working with the difference schemes. Recall from Section~\ref{basicdif} that a difference scheme is {\em
stable} if the corresponding difference operators $R_h$ (that send
the grid function $[u^l]$ to the grid
function $[u^{l+1}]$) are bounded
uniformly on $h$, together with  their powers. The investigation
of stability of the difference scheme from the previous subsection
can be done as follows (for more details see e.g.\ \cite[p.\ 78]{go76}).

First  consider the one-dimensional scheme~\eqref{scheme_1} and
the case $\mu\geq 0$ (in case $\mu<0$ the argument is similar). Denote
by $\nu=|\mu|\displaystyle\frac{\tau}{h}$ the Courant number and
check the scheme stability by the Fourier method. Substituting in~
\eqref{scheme_1} the values (where $i^2=-1$)
 $$w_{j-\frac{1}{2}}=w^*e^{ij\phi},\
\ w^{j-\frac{1}{2}}=\lambda w_{j-\frac{1}{2}},$$
 we obtain the characteristic equation
 $$\lambda(\phi)=1-\nu(1-e^{-i\phi}).$$
 The necessary and sufficient condition for the
one-dimensional difference operator to be uniformly bounded,
together with its powers,  is the condition $|\lambda(\phi)|\leq
1$ for all $\phi\in[0,2\pi),$ that is equivalent to the condition
$\nu\leq 1$ for the Courant number (it follows from a rather technical, but simple argument).

For the $n$-dimensional scheme~\eqref{scheme_2_sist} approximating
the boundary-value problem~\eqref{sist_2} (when $m=2$), the
stability condition is
 \begin{equation}\label{ustoj}
 \tau\left(\frac{1}{\tau_x}+\frac{1}{\tau_y}\right)\leq 1,
 \end{equation}
where $\tau_x,
\tau_y$ as in~\eqref{stab} are the
maximal time steps guaranteeing the stability of the corresponding
one-dimensional schemes.

The proof of this stability condition is based on the case of a
one-dimensional scheme  for one equation (described above), on the theory of matrix pencils $\lambda A-B_i$ \cite{gant,hj} and on
a difference energy integral inequality: under the restriction~
\eqref{ustoj} one has
$$\sum\limits_{i,j=1}^{2^N}(A{\bf
u}^{i-\frac{1}{2},j-\frac{1}{2}},{\bf
u}^{i-\frac{1}{2},j-\frac{1}{2}})-\sum\limits_{i,j=1}^{2^N}(A{\bf
u}_{i-\frac{1}{2},j-\frac{1}{2}},{\bf
u}_{i-\frac{1}{2},j-\frac{1}{2}})\leq$$
\begin{equation}\label{energy1}\frac{\tau}{h}\left(\sum\limits_{j=1}^{2^N}\left[
(B{\mathcal U}_{0,j-\frac{1}{2}},{\mathcal U}_{0,j-\frac{1}{2}})-
(B{\mathcal U}_{2^N,j-\frac{1}{2}},{\mathcal U}_{2^N,j-\frac{1}{2}})\right]\right)+\end{equation}
$$
+\frac{\tau}{h}\left(\sum\limits_{i=1}^{2^N}\left[ (C{\mathcal
U}_{i-\frac{1}{2},0},{\mathcal U}_{i-\frac{1}{2},0})- (C{\mathcal
U}_{i-\frac{1}{2},2^N},{\mathcal
U}_{i-\frac{1}{2},2^N})\right]\right).$$
 
Due to the dissipativity of the boundary conditions, the
right-hand part of the inequality~\eqref{energy1} is below zero, thus
\begin{equation}\label{ust}\sum\limits_{i,j=1}^{2^N}(A{\bf
u}^{i-\frac{1}{2},j-\frac{1}{2}},{\bf
u}^{i-\frac{1}{2},j-\frac{1}{2}})\leq\sum\limits_{i,j=1}^{2^N}(A{\bf
u}_{i-\frac{1}{2},j-\frac{1}{2}},{\bf
u}_{i-\frac{1}{2},j-\frac{1}{2}})\end{equation}
 which is equivalent to the
stability condition.

We omit the proof of the energy inequality~\eqref{energy1} since it is standard (see e.g.\ \cite{go76,fri,gr,gv}) and rather technical. \qedhere
\end{enumerate}
\end{proof}

\section{Computability of the solution
operators}\label{smain}

In this section we give precise formulations  and proof schemes of our
main results. The precise formulations are given in Subsection~\ref{rmain}.
In Subsections~\ref{scmain} and~\ref{scmain2} we describe the proof
schemes of the main results omitting technical details of the
relevant estimates. The technical details are presented in the
next section.

\subsection{Formulations of main results}\label{rmain}

Let us formulate the main results of this paper.  The first main result
concerns  computability of the  boundary-value problem~
\eqref{sist_2} posed in Section~\ref{state}.

\begin{thm}\label{main}
Let $Q=[0,1]^m$; $T>0$ be a computable real and $M_{\varphi}>0$, $p\geq 2$ be integers.
Let  $A, B_1,\ldots,B_m$ be fixed computable symmetric
matrices, such that $A=A^\ast> 0$, $B_i=B_i^\ast$.
 Let
$\Phi_i^{(1)}$, $\Phi_i^{(2)}$ $(i=1,2,\ldots,m)$ be fixed
computable  rectangular real non-degenerate
matrices, with their numbers of rows equal to the number of
positive and negative
 eigenvalues of $A^{-1}B_i$,
respectively, and such that inequalities~\eqref{gran_dissip} hold.

If $\varphi\in C^{p+1}(Q)$ satisfies
 \begin{equation}\label{bound}
||\frac{\partial
\varphi}{\partial x_i}||_{s}\leq M_{\varphi},\ ||\frac{\partial^2
\varphi}{\partial x_i\partial x_j}||_{s}\leq M_{\varphi},\
 i,j=1,2,\ldots,m,
 \end{equation}
 and meets the boundary conditions,  then the operator
$R:\varphi \mapsto {\bf u}$ mapping the initial function 
to the unique solution ${\mathbf
u}\in C^p(Q\times[0,T],{\mathbb R}^n)$ of the
boundary-value problem~\eqref{sist_2} is a computable partial
function from $C_s(Q,\mathbb R^n)$ to $C_{sL_2}(Q\times[0,T],\mathbb
R^n)$.
\end{thm}

A natural question is whether the computability of solution operator $R$ is
uniform  on the matrices $A,B_i$. Currently we do not know the answer for the
arbitrary real matrices  (cf.\ Remark~\ref{comcoef}.2 below) but the uniformity holds when the coefficients of $A,B_i$ range through an arbitrary strongly  constructive real
closed ordered subfield $(\mathbb{B},\beta)$ of $\mathbb{R}_c$. The next result extends the previous theorem because, by Lemma~\ref{rfields}, for any computable real matrices $A,B_i,\Phi_i^{(1)}$, $\Phi_i^{(2)}$ there is a strongly  constructive real
closed ordered subfield $(\mathbb{B},\beta)$ of $\mathbb{R}_c$
that contains all coefficients of the matrices $A,
B_1,\ldots,B_m,\Phi_i^{(1)}$, $\Phi_i^{(2)}$. 

\begin{thm}\label{commain}
Let $(\mathbb{B},\beta)$ be a strongly  constructive real
closed ordered subfield  of $\mathbb{R}_c$. Then the operator $R$ from the previous theorem is uniformly computable (w.r.t the numbering $\beta$) on the matrices $A,
B_1,\ldots,B_m,\Phi_i^{(1)}$, $\Phi_i^{(2)}$ with coefficients in $\mathbb{B}$. 
 \end{thm}

The second main result concerns the initial value problem~
\eqref{sist_1}. It improves the main result of \cite{ss09} and is uniform on the arbitrary matrices $A,B_i$.

\begin{thm}\label{main2}
Let $M_{\varphi}>0, M_A>0, p\geq2$ be integers,   let
$i=1,\ldots,m$, and let $n_A,n_1,\ldots,n_m$ be cardinalities of spectra of $A$ and of the matrix pencils $\lambda A- B_1,\ldots,\lambda A-B_m$, respectively (i.e., $n_i$ is the number of distinct roots of the characteristic polynomial
${\rm det}(\lambda A- B_i)$). Then the operator
 $$(A,B_1,\ldots,B_m,n_A,n_1,\ldots,n_m,\varphi)\mapsto{\mathbf u}$$
  sending any  sequence
$A,B_1,\ldots,B_m$ of symmetric real matrices with $A> 0$ such that the matrix pencils $\lambda A- B_i$ have no zero eigenvalues,
 \begin{equation}\label{M_A}
 ||A||_2,||A^{-1}||_2,||B_i||_2\leq M_A,\quad
\lambda^{(i)}_{\min}<0<\lambda^{(i)}_{\max},\quad i=1,2,\ldots,m,
\end{equation}
the sequence $n_A, n_1, \ldots, n_m$ of the corresponding cardinalities, and any
 $\varphi \in C^{p+1}(Q,{\mathbb R}^n)$ satisfying
the conditions~\eqref{bound}, to the unique solution ${\mathbf
u}\in C^p(H,{\mathbb R}^n)$ of~\eqref{sist_1} is a computable
partial function from the space $S_{+}\times S^{m}\times \mathbb{N}^{m+1}\times
C_s(Q,{\mathbb R}^n)$ to $C_{sL_2}(H,{\mathbb R}^n)$.
\end{thm}

In Theorem~\ref{main2},  $\lambda^{(i)}_{\min},
\lambda^{(i)}_{\max}$  are respectively the minimal and maximal
eigenvalues of the matrix pencil $\lambda A- B_i$,  $S\subseteq
\mathbb{R}^{n\times n}$ is the space of symmetric $n\times n$
matrices equipped with the Euclidean norm, $S_{+}$ is the space of
symmetric positively definite matrices with the Euclidean norm,
and $H\subseteq{\mathbb R}^{m+1}$ is the domain of correctness of~
\eqref{sist_1}, i.e., the maximal set where, for any $p\geq 2$ and
$\varphi\in C^{p+1}(Q,{\mathbb R}^n)$, there exists a unique
solution ${\bf u}\in C^p(H,{\mathbb R}^n)$ of the initial value
problem~\eqref{sist_1}.

The set $H$ is known to be (see e.g.\ 
\cite{go71}) a nonempty intersection of the semi-spaces
 $$
t\geq0,\;x_i-\lambda^{(i)}_{\rm
max}t\geq0,\;x_i-1-\lambda^{(i)}_{\rm min}t\leq0,\;(i=1,\ldots,m)
 $$
of ${\mathbb R}^{m+1}$. We are especially interested in the case when $H$ is a
compact subset of $Q\times~[0,+\infty)$ (obviously, a sufficient
condition for this to be true is $\lambda^{(i)}_{\rm
min}<0<\lambda^{(i)}_{\rm max}$ for all $i=1,\ldots,m$; this is
often the case for natural physical systems.
In \cite{ss09} we observed that the domain $H$ for the problem~\eqref{sist_1}
is computable from $A,B_1,\ldots,B_m$ (more exactly, the vector
$(\lambda^{(1)}_{\max},\ldots,\lambda^{(m)}_{\max},
\lambda^{(1)}_{\min},\ldots,\lambda^{(m)}_{\min})$  is computable
from $A,B_1,\ldots,B_m$; this implies  computability of $H$ in
the sense of computable analysis \cite{wei}).

Since, for each $i=1,\ldots,m$, $\lambda^{(i)}_{\max}$ is the
maximal and $\lambda^{(i)}_{\min}$ is the minimal eigenvalue of
the matrix pencil $\lambda A-B_i$, and maximum and minimum of a
vector of reals are computable \cite{wei}, it suffices to show
that a vector $(\lambda_1,\ldots,\lambda_n)$  consisting  of all
eigenvalues of a matrix pencil $\lambda A-B$ is computable from
$A,B$. But $(\lambda_1,\ldots,\lambda_n)$ is a vector of all roots
of the characteristic polynomial of $\lambda A-B$, hence it is
computable \cite{wei,bhw}. We immediately obtain

\begin{lem}\label{main11}
A rational  number $\tau$  meeting~\eqref{ustoj} is
computable from symmetric real matrices $A,B_1,\ldots,B_m$ such
that $A>0$.
\end{lem}

\begin{rems}\label{comcoef}\leavevmode
  \begin{enumerate}[label=\arabic*.]
    \item Besides  the condition
      $\lambda^{(i)}_{\min}<0<\lambda^{(i)}_{\max}$ in Theorem~ \ref{main2},
      some  alternative natural conditions may be assumed.
E.g., for one equation $\frac{\partial u}{\partial t}-\frac{\partial
u}{\partial x}=0$ the domain of correctness may be the intersection of the semi-planes
$\{t\geq 0\}$, $\{x\leq t\}$, $\{x\leq 1+t\}$, and we search the
solution in the intersection of the semi-planes $\{t\geq 0\}$, $\{x\leq t\}$, $\{x\leq
1\}$. Our proof is adjusted to similar modifications in a
straightforward way.

\item Note that  Theorem~\ref{main2} states   computability on
arbitrary matrices $A,B_i$ while Theorem~\ref{main} (and Theorem~\ref{commain}) does not. The reason
is that our proof of Theorem~\ref{main2} (where we take rational
fast Cauchy approximations to $A,B_i$) can not be
straightforwardly adjusted to that of Theorem~\ref{main} because
the dissipativity conditions in the last theorem might hold for
the given matrix but not hold for the approximate matrices.
Currently we do not know whether Theorem~\ref{commain} may be
strengthened to include  computability on the real coefficients
$A,B_i$.
If we strengthen the dissipativity conditions to strict inequalities then the
solution operator in Theorem~\ref{main} will be computable on
$A,B_i,\Phi_i^{(1)}$, $\Phi_i^{(2)}$ $(i=1,2,\ldots,m)$, similarly to Theorem~\ref{main2}. 

\item In \cite{ss09}, we established a weaker result that the solution is computable provided that $A,B_1,\ldots,B_m$ are fixed computable matrices (in this case one can of course omit the conditions on spectra of $A$ and of the matrix pencils). This weaker result is proved just in the same way as Theorem~\ref{main} below. In the stronger formulation above, the proof requires additional considerations described below.
\end{enumerate}

 \end{rems}

\subsection{Scheme of  proof of Theorem~\ref{main}}\label{scmain}

Here we provide outline of the proof of Theorem~\ref{main} (which also applies to Theorem~\ref{commain})
omitting the proofs of some technical estimates. The estimates are
proved in the next section.

\begin{enumerate}[label=\arabic*.]
\item Consider a grid $G=G_k$ of step $h=2^{-k}$ on $Q$, as described
in Subsection~\ref{spaces}, and choose a computable sequence
$\{T_k\}$ of rational numbers that fast  converges to $T$.  Take a
sequence $\{\varphi_k\}$ of grid functions $\varphi_k:
G_{i_k}\rightarrow{\mathbb Q}^n$  such that their multilinear
interpolations $\{\tilde{\varphi}_k\}$ form a fast Cauchy sequence
converging in $C_s(Q,{\mathbb R}^n)$ to $\varphi$, so
 \begin{equation}
||\tilde{\varphi}_k-\varphi||_s\leq\frac{1}{2^k.}
 \end{equation}

\item Choose $(\mathbb{B},\beta)$ as in Remark~\ref{commain}. We compute a
sequence $\{\upsilon_k\}$ of grid functions
${\mathbf\upsilon}_k:G^\tau_{i_k}\rightarrow{\mathbb B}^n$ (for some sequence
$\{i_k\}$ of natural numbers). Without loss of generality we may
assume that the sequence $\{i_k\}$ is increasing (otherwise, choose
a suitable subsequence of $\{\varphi_k\}$). Let the grid function
$\upsilon_k$ be constructed from $\varphi_k$, $A$, $B_1,\ldots,B_m$
and $\Phi_i^{(1)}$,  $\Phi_i^{(2)}$ $(i=1,2,\ldots,m)$, by the
algorithm of the difference equation in Subsection~\ref{condif}.
According to Lemma~\ref{complin}, the operation $\varphi_k\mapsto
\upsilon_k$ is computable (w.r.t. $\beta$) and linear on
$\varphi_k$.

 By Lemma~\ref{fspace}, it
suffices to show that for some constant $c$ (depending only on
$A,B_1,\cdots,B_m$, $\Phi_i^{(1)}$, $\Phi_i^{(2)}$ and  $M_{\varphi}$)
we have
\begin{equation}\label{estim}
 ||\tilde{\upsilon}_k-{\bf u}||_{sL_2}\leq
c\cdot\frac{1}{2^k}.
 \end{equation}
 for all $k$,
i.e.\ $\{\tilde{\upsilon}_k\}$ fast converges in
$C_{sL_2}(H,{\mathbb R}^n)$ to ${\bf u}$.

\item We divide the proof of~\eqref{estim} into several parts. For
any $k$, let $\tilde{u}_k$ be the interpolation  of  the grid
function computed by the algorithm in Subsection~\ref{difop} from
the exact initial values $\varphi|_{G_k}$. By $\widetilde{
u|_{G_k^{\tau_k}}}$ we denote the interpolation of the
$G_k$-discretization of the solution ${\bf u}$ of the differential
problem~\eqref{sist_2}. We will estimate independently the
following three summands:
\begin{equation}\label{estim_three}
||\tilde{\upsilon}_k-{\bf u}||_{sL_2}\leq
||\tilde{\upsilon}_k-\tilde{u}_k||_{sL_2}+ ||\tilde{u}_k-\widetilde{
u|_{G_k^{\tau_k}}}||_{sL_2}+||\widetilde{
u|_{G_k^{\tau_k}}}-{\bf u}||_{sL_2}.
 \end{equation}

\item The third summand is estimated with the help of the properties
of interpolations:
 \begin{equation}\label{3}||\widetilde{
u|_{G_k^{\tau_k}}}-{\bf u}||_{sL_2}\leq
c_{int}\cdot 2^{-k}.
 \end{equation}
 The constant $c_{int}$ depends on the norms
of derivatives of ${\bf u}$, which can be  estimated by $\varphi$
and its derivatives (hence by some expression involving
$M_{\varphi}$ and the norms of the coefficient matrices, following
the lines of the proof of the uniqueness theorem for~\eqref{sist_2}
in \cite[p.\ 194]{go76}, see also \cite{evans}).

\item The second summand is estimated with the help of Theorem~\ref{convds}
  on convergence of the difference scheme in grid norms:
 \begin{equation}\label{2}
 ||\tilde{u}_k-\widetilde{
u|_{G_k^{\tau_k}}}||_{sL_2}\leq c_{diff}\cdot 2^{-k}.
 \end{equation}
 The constant $c_\mathit{diff}$ also depends only on the
derivatives of $\varphi$ and on the coefficients of~\eqref{sist_2}.

\item The first summand  is estimated by means of  stability of
the difference scheme, with taking into account linearity of the
difference and interpolation operators. The corresponding constant
depends on the coefficients of the differential equations:
 \begin{equation}\label{1}
 ||\tilde{\upsilon}_k-\tilde{u}_k||_{sL_2}\leq c_A\cdot 2^{-k}
 \end{equation}
\end{enumerate}

\subsection{Scheme of  proof of Theorem~\ref{main2}}\label{scmain2}

Here we provide outline of the proof of Theorem~\ref{main2}
omitting the proofs of some technical estimates. The estimates are
proved in the next section.

\begin{enumerate}[label=\arabic*.]
\item  Consider again the grid $G=G_k$ of step $h=2^{-k}$ on $Q$ .Take
a sequence $\{\varphi_k\}$ of grid functions $\varphi_k:
G_{i_k}\rightarrow{\mathbb Q}^n$  such that their multilinear
interpolations $\{\tilde{\varphi}_k\}$ form a fast Cauchy sequence
converging in $C_s(Q,{\mathbb R}^n)$ to $\varphi$.

\item Let $A^{(k)}$ and $B_i^{(k)}$  be sequences of symmetric
matrices, such that $A^{(k)}>0$, with rational coefficients fast
converging to $A$ and $B_i$ respectively, $i=1,2,\ldots,m$, in the standard
Euclidean norm 
 \begin{equation}
 ||A-A^{(k)}||_{2}\leq 2^{-k},\;||B_i-B_i^{(k)}||_{2}\leq 2^{-k}.
 \end{equation}

From \cite[Theorem~2]{zb04} (stating that, given a complex normal matrix
and the cardinality of its spectrum one can compute the sequence of
all its eigenvalues counted with their multiplicities, as well as an orthonormal basis of eigenvectors), we can
without loss of generality assume that, for all $k$ and $i=0,\ldots, m$,
the matrix pencil $\lambda A^{(k)}-B^{(k)}_i$ has no zero
eigenvalues.
By item~(2) of Lemma~\ref{complin}, the function
$(A,B_1,\ldots,B_m,k)\mapsto\tau$ is computable, so there is a
computable sequence $\{\tau_k\}$ of positive rationals that fast
converges to $\tau$ and, for any $k$, $\tau_k$ satisfies the
stability condition~\eqref{ustoj} for matrices $A,B_1,\ldots,B_m$,
and also for matrices $A^{(k)},B^{(k)}_1,\ldots,B^{(k)}_m$.

Let $\{\upsilon_k\}$ be an $\alpha$-computable sequence obtained as in the proof
of Theorem~\ref{main},
only with the strongly constructive ordered field
$(\mathbb{A},\alpha)$ of algebraic reals in place of
$(\mathbb{B},\beta)$ and with $A^{(k)}$, $B_i^{(k)}$ in place of
$A,B_i$. It suffices to show that for some constant $c$ (depending
only on $M_A$ and  $M_{\varphi}$) we have~\eqref{estim}
 for all $k$,
i.e.\ $\{\tilde{\upsilon}_k\}$ fast converges in
$C_{sL_2}(H,{\mathbb R}^n)$ to ${\bf u}$.

\item We divide the proof of~\eqref{estim} again into several parts.

 Let
$\widehat{\upsilon}_k$ be the grid function obtained by the
difference scheme of Subsection~\ref{condif} from $\varphi_k$ and
$A,B_i$ (the ``exact'' coefficient matrices); let
$\widetilde{\widehat{\upsilon}_k}$ be its interpolation.

Let $\tilde{u}_k$ be the sequence of  interpolations of the grid
functions obtained by the difference scheme of Subsection~\ref{condif}
from the exact initial values $\varphi|_{G_k}$  and
from the exact matrices $A,B_i$. By $\widetilde{
u}|_{G_k^{\tau_k}}$ we  denote the interpolation of the
$G_k$-discretization of the solution ${\bf u}$ of the differential
problem~\eqref{sist_1}. Now~\eqref{estim} is naturally splitted to
four summands:
\begin{equation}\label{estim_four}
||\tilde{\upsilon}_k-{\bf u}||_{sL_2}\leq ||\tilde{\upsilon}_k-\tilde{\widehat{\upsilon}}_k||_{sL_2}+
||\tilde{\widehat{\upsilon}}_k-\widetilde{u}_k||_{sL_2}+ ||\tilde{u}_k-\widetilde{
u}|_{G_k^{\tau_k}}||_{sL_2}+||\widetilde{
u}|_{G_k^{\tau_k}}-{\bf u}||_{sL_2}.
 \end{equation}

\item The last three summands are estimated in the same way as~
\eqref{3}--\eqref{1} in the previous subsection.  The only
difference is that all estimates are obtained in the domain $H$ of
uniqueness of the Cauchy problem rather than in the set $Q\times
[0,T]$ in the proof of Theorem~\ref{main}. The procedure of
identifying the grid cells, which are in $H$, was described in
\cite{ss09}.

\item The key technical tool for estimating the first summand
 \begin{equation}\label{11}
||\tilde{\upsilon}_k-\tilde{\widehat{\upsilon}}_k||_{sL_2}\leq c_{rat}\cdot
2^{-k}
 \end{equation}
 is formal differentiation of the difference scheme. For doing this correctly,
the assumption of Theorem~\ref{main2} that the eigenvalues are
non-zero  and the cardinalities of spectra are known in advance, are
needed. Note that the constant $c_{rat}$ depends on the eigenvalues
of the matrices $A, B_i$.
\end{enumerate}

\section{Proofs of the estimates}\label{estimates}

In this section we prove the technical estimates from the previous
section.

\subsection{Interpolation and proof of the estimates~\eqref{3},~\eqref{1}}\label{estim3}

Recall the construction of the multilinear interpolations.

In the one-dimensional case, the interpolating function
$\tilde{\bf u}$ is defined inside the grid rectangles
 $$
\left(i-\frac{1}{2}\right)h\leq
x\leq\left(i+\frac{1}{2}\right)h;\;l\tau\leq t\leq(l+1)\tau
 $$
in the standard way as follows
  \begin{eqnarray*}
\tilde{\bf u}(x,t)&=&{\bf
u}_{i-\frac{1}{2}}\cdot\left(l+1-\frac{t}{\tau}\right)\cdot\left(i+\frac{1}{2}-\frac{x}{h}\right)
 +{\bf  u}_{i+\frac{1}{2}}\cdot\left(l+1-\frac{t}{\tau}\right)\cdot\left(\frac{x}{h}-\left(i-\frac{1}{2}\right)\right)\\
 &+&{\bf u}^{i-\frac{1}{2}}\cdot\left(\frac{t}{\tau}-l\right)\cdot\left(i+\frac{1}{2}-\frac{x}{h}\right)
 +{\bf u}^{i+\frac{1}{2}}\cdot\left(\frac{t}{\tau}-l\right)\cdot\left(\frac{x}{h}-\left(i-\frac{1}{2}\right)\right)
  \end{eqnarray*}
where ${\bf u}_{i\pm\frac{1}{2}}$ and ${\bf u}^{i\pm\frac{1}{2}}$
are the grid functions on time levels $t=l\tau$ and $t=(l+1)\tau$,
respectively.

In the two-dimensional case (and for higher dimensions $m$) the
interpolating function is defined in a similar way. Since the full
expression is rather long we write down only two (of eight)
summands, others are constructed in an obvious way (see \cite[p.\ 212]{go71}):
  \begin{eqnarray*}
\tilde{\bf u}(x,y,t)&=&{\bf
u}_{i-\frac{1}{2},j-\frac{1}{2}}\cdot\left(l+1-\frac{t}{\tau}\right)\cdot
\left(i+\frac{1}{2}-\frac{x}{h}\right)\cdot\left(j+\frac{1}{2}-\frac{y}{h}\right)\\
 &+&{\bf
 u}_{i+\frac{1}{2},j-\frac{1}{2}}\cdot\left(l+1-\frac{t}{\tau}\right)\cdot
 \left(\frac{x}{h}- \left(i-\frac{1}{2}\right)\right)\cdot\left(j+\frac{1}{2}-\frac{y}{h}\right)\\
 &+&\cdots
  \end{eqnarray*}
where $\left(j-\frac{1}{2}\right)h\leq
y\leq\left(j+\frac{1}{2}\right)h$.

From these formulas for multilinear interpolation, linearity of the
interpolation operators ${\bf u}\mapsto\tilde{{\bf u}}$ and
${\bf u}|_{G_k}\mapsto\widetilde{{\bf u}|_{G_k}}$ is obvious.

\begin{prop}\label{tildeu}
$||\widetilde{u|_{G_k}}-{\bf u}||_{sL_2}\leq
c_{int}\cdot\frac{1}{2^k}$ for some constant $c_\mathit{int}$ depending only on $M_A,M_\varphi$.
 \end{prop}

\begin{proof}  By a well-known estimate for the multilinear
interpolations \cite{sz59,so74,za,ba86},
\begin{equation}\label{est-int}||\widetilde{u|_{G_k}}-{\bf u}||_s\leq c_{int}\cdot\frac{1}{2^k}\end{equation} for
some constant $c_{int}$ depending only on the $s$-norms of second
derivatives of ${\bf u}$. Since the $s$-norm is stronger than the
$sL_2$-norm, $||\tilde{u}_k-{\bf u}||_{sL_2}\leq
c_{int}\cdot\frac{1}{2^k}$. It suffices to show that $c_{int}$
depends in fact only on the  second derivatives of $\varphi$ and the
norms of the matrices $A$, $A^{-1}$, $B_i$.

Considering the Cauchy problem,  due to
the smoothness assumptions, we can construct auxiliary Cauchy
problems for partial derivatives of ${\bf u}$ (we write down a couple of
them, as examples):
 $$
\begin{cases} A({\bf u}_x)_t+B({\bf u}_x)_x+C({\bf u}_x)_y=0,\\
{\bf u}_x|_{t=0}=\varphi_x,
\end{cases}
 $$
 $$
\begin{cases} A({\bf u}_t)_t+B({\bf u}_t)_x+C({\bf u}_t)_y=0,\\
{\bf u}_t|_{t=0}=-A^{-1}(B\varphi_x+C\varphi_y),
\end{cases}
 $$
 \begin{equation} \label{sist_13}
\begin{cases} A({\bf u}_{tt})_t+B({\bf u}_{tt})_x+C({\bf u}_{tt})_y=0,\\
{\bf u}_{tt}|_{t=0}=-A^{-1}(B({\bf u}_t|_{t=0})_x+C({\bf u}_t|_{t=0})_y)=\psi(x,y).
\end{cases}
\end{equation}
 As it is known, we have 
 \begin{equation}\label{energy}||{\bf u}||_{A,sL_2}\leq ||\varphi||_{A,L_2}\end{equation} 
 The proof of this estimate is rather long and technical; it is presented in
 detail for the considered systems of PDEs (even with a nonzero righthand part
 $f$; in this case the estimate contains also the norm of $f$) in \cite{go71}
 (the estimates are stated on p.\ 155 for the Cauchy problem and on p.\ 194 for
 the boundary-value problem, respectively), and also can be found in
 \cite[subsection~7.3]{evans}. Applying an analog of~\eqref{energy}
 to the systems for the second derivatives of {\bf u} and using the equivalence
 of norms in $\mathbb R^n$, we obtain
 \begin{align*}
   ||\frac{\partial^2{\bf u}}{\partial x_i\partial x_j}||_{sL_2}\leq
\sqrt{\frac{\lambda_{max}(A)}{\lambda_{min}(A)}}\cdot
c(||A||_2,||A^{-1}||_2,||B_1||_2,\ldots,||B_1||_2,\\
{\rm max}_{i,j=1,\ldots,m}||\frac{\partial^2\varphi}{\partial
x_i\partial x_j}||_{L_2}) \leq c(M_A,M_\varphi). \tag*{\qEd}
 \end{align*}
\def\popQED{}
 \end{proof}

Thus the estimate~\eqref{3} needed in the proofs of Theorems 2 and 3 is established.
Further we will also need the following

\begin{lem}\label{tildeuh}
Let $u^{(h)}$ be calculated from $\varphi^{(h)}$ by means of the
difference scheme in Subsection~\ref{condif}. Then
$||\widetilde{u^{(h)}}||_{sL_2}\leq
\sqrt{\frac{\lambda_{max}(A)}{\lambda_{min}(A)}}||\widetilde{\varphi^{(h)}}||_s$.
 \end{lem}

\proof Let us first show that
 \begin{equation}\label{14}
\max_{0\leq l\tau\leq
T}\int\limits_{H\cap\{t=l\tau\}}|\tilde{u}(x,y,t)|^2dxdy\leq
\max_{0\leq l\tau\leq T}\left(h^2\sum\limits_{i,j} u^2_
{i-\frac{1}{2},j-\frac{1}{2}}\right)=||u^{(h)}||^2_{sL_2}.
 \end{equation}
 For simplicity of notation consider the one-dimensional case
(adding an additional variable is straightforward). In the $i$-th
grid cell for a fixed $t=l\tau$ we have
 $$\tilde{u}(x,t)=u_{il}\cdot(i+1-\frac{x}{h})+u_{il+1}\cdot(\frac{x}{h}-i)
 $$
 and
 \begin{eqnarray*}\label{140}
\int\limits^{(i+1)h}_{ih}\tilde{u}^2(x,t)dx=
\int\limits^{(i+1)h}_{ih}[u_{j,l}(i+1-\frac{x}{h})+u_{i,l+1}(\frac{x}{h}-i)]^2dx=\\
h\int\limits^1_0[u_{i,l}(1-\xi)+u_{i,l+1}\xi]^2d\xi=\frac{h}{3}(u_{i,l}^2+
 u_{i,l}u_{j,l+1}+u_{i,l+1}^2)\leq h\frac{u_{i,l}^2+
 u_{i,l+1}^2}{2}.
 \end{eqnarray*}
 The summation by $i$ yields
 $\int\limits^1_0\tilde{u}^2(x,l\tau)dx\leq h\sum\limits_iu^2_{i,l}(x,l\tau)$.
Taking maximum over all $l$  concludes the proof
of~\eqref{14}.

It is well-known \cite{sz59,so74}, that  for the linear interpolations $||\widetilde{u^{(h)}}||_{sL_2}\leq||u^{(h)}||_{sL_2}$ where the
right-hand part refers to the grid norm. Obviously,
 \begin{equation}\label{141}
||\varphi^{(h)}||^2_{L_2}=h^2\sum\limits_{i,j}\varphi^2_{i-\frac{1}{2},j-\frac{1}{2}}\leq
h^2\frac{1}{h^2}\max_{i,j}\varphi^2_{i-\frac{1}{2},j-\frac{1}{2}}\leq
\sup_{(x,y)\in Q}\widetilde{\varphi|_G}^2(x,y)=||\varphi^{(h)}||^2_s.
 \end{equation}
 The estimate~\eqref{ust} implies
$||u^{(h)}||_{A,sL_2}\leq||\varphi^{(h)}||_{A,L_2}$. Taking into
account the equivalence of the Euclidean norms
$\lambda_{min}(A)\langle {\bf u},{\bf u}\rangle\leq\langle
A{\bf u},{\bf u}\rangle\leq\lambda_{max}(A)\langle {\bf u},{\bf u}\rangle$ we obtain the
desired estimate.
\qed

Arguments similar to those in the proof of Lemma~\ref{tildeuh} can be found
in \cite{go71}. We have recalled them for the convenience of the
reader.

\begin{prop} The estimate~\eqref{1} holds.
 \end{prop}

\proof Using linearity of the interpolation and difference
operators and Lemma~\ref{tildeuh} we obtain
 \begin{eqnarray*}\label{ocenka11}
||\tilde{\upsilon}_k-\tilde{u}_k||_{sL_2}=||\widetilde{\upsilon_k-u_k}
||_{sL_2}
\leq
\sqrt{\frac{\lambda_{max}(A)}{\lambda_{min}(A)}}||\widetilde{\varphi_k} -\widetilde{\varphi|_{G_k}}||_{s}\\
\leq
\sqrt{\frac{\lambda_{max}(A)}{\lambda_{min}(A)}} (||\widetilde{\varphi_k}-\varphi||_s+
||\widetilde{\varphi|_{G_k}}-\varphi||_s)\leq
c(M_A,M_{\varphi})2^{-k}.
 \end{eqnarray*}
 Here, $||\widetilde{\varphi|_{G_k}}-\varphi||_s$ is again estimated
 by~\eqref{est-int}.
\qed

\subsection{Convergence of the difference scheme and proof of the estimate~\eqref{2}}\label{estim3}

\begin{lem}\label{conv}
There is a constant $c_{\operatorname{diff}}$ depending only on
$M_A,M_\varphi$ such that for all $k\geq0$ we
have $||u_k-u|_{G_k^\tau}||_{sL_2}\leq
c_{\mathit{diff}}\cdot\frac{1}{2^k}$ where $u$ is the
solution of~\eqref{sist_1} or~\eqref{sist_2}.
\end{lem}

\proof  The  estimate follows from Theorem~\ref{convds}. The
fact that $c_{\mathit{diff}}$ depends only on $M_A$ and
$M_{\varphi}$ follows from the  proof of this theorem in
\cite[Chapter 5]{gr} according to which we can take
$c_{\operatorname{diff}}=c_1\cdot c_2$ where $c_2$ comes from the
stability condition and $c_1$ is from the approximation
$||L_hu_h-(Lu)|_{G^\tau_k}||_{sL_2}\leq c_1h$. Since we consider a
first-order difference scheme, it follows from the Taylor
decomposition of ${\bf u}$ that $c_1$ depends only on $M_A$ and
 $$||\frac{\partial {\bf u}}{\partial x_i}||_{sL_2},\; ||\frac{\partial
{\bf u}}{\partial t}||_{sL_2},\; ||\frac{\partial^2 {\bf u}}{\partial
x_i\partial x_j}||_{sL_2},\; ||\frac{\partial^2 {\bf u}}{\partial
x_i\partial t}||_{sL_2}.$$
 By the proof of the uniqueness theorem for~
\eqref{sist_1} (see the
proof of Proposition~\ref{tildeu}), the norms of the derivatives
above are bounded by a constant depending only on $M_A,M_\varphi$.
 \qed

\begin{prop}\label{conv1}
There is a constant $c$ depending only on
$M_A,M_\varphi$ such that for all $k\geq0$ we
have
$||\tilde{u}_k-\widetilde{u|_{G_k}}||_{sL_2}\leq c\cdot\frac{1}{2^k}$ where $u$ is the
solution of~\eqref{sist_1} or~\eqref{sist_2}.
\end{prop}

\proof  Since the  operator of multilinear interpolation is linear, from~
\eqref{14} and the estimate of the previous lemma we obtain
 $$||\tilde{u}_k-\widetilde{u|_{G_k^\tau}}||_{sL_2}\leq
||u_k-u|_{G_k^\tau}||_{sL_2}\leq c_{diff}2^{-k}.
 $$
 This implies the desired estimate,  which is exactly the
estimate~\eqref{2} needed in the proofs of Theorems~\ref{main},~\ref{main2}.
 \qed

\subsection{Proof of the  estimate~\eqref{11}}\label{estim4}

Finally, we prove the  estimate~\eqref{11} which is needed only  for
the Cauchy problem. Recall that $\upsilon_k$ and $\widehat{\upsilon}_k$
satisfy respectively the following difference schemes (see~
\eqref{scheme_2_sist}) in which the index $k$ in $\tau_k$,
$\upsilon_k$, $\widehat{\upsilon}_k$ and $w_k$ is omitted  for
simplicity (note that $\tau_k$ is chosen in such a way that both
difference schemes below are stable)
 \begin{equation}\label{sist_dif}
 A\frac{\widehat{{\bf \upsilon}}^{i-\frac{1}{2},j-\frac{1}{2}}-\widehat{{\bf
\upsilon}}_{i-\frac{1}{2},j-\frac{1}{2}}}{\tau}+
B\frac{\widehat{\Upsilon}_{i,j-\frac{1}{2}}-\widehat{
\Upsilon}_{i-1,j-\frac{1}{2}}}{h}+C\frac{\widehat{
\Upsilon}_{i-\frac{1}{2},j}-\widehat{
\Upsilon}_{i-\frac{1}{2},j-1}}{h}=0,
 \end{equation}
 \begin{equation}\label{sist_dif_hat}
A^{(k)}\frac{{\bf \upsilon}^{i-\frac{1}{2},j-\frac{1}{2}}-{\bf
\upsilon}_{i-\frac{1}{2},j-\frac{1}{2}}}{\tau}+B^{(k)}\frac{
\Upsilon_{i,j-\frac{1}{2}}-
\Upsilon_{i-1,j-\frac{1}{2}}}{h}+C^{(k)}\frac{
\Upsilon_{i-\frac{1}{2},j}-
\Upsilon_{i-\frac{1}{2},j-1}}{h}=0
 \end{equation}
 with the initial conditions  $\upsilon|_{t=0}=\widehat{{\bf \upsilon}}|_{t=0}=\varphi_k$.
Deducting the second system of  equations from the first one  we obtain the difference equations
\begin{equation}\label{tilde}
 A\frac{{\bf w}^{i-\frac{1}{2},j-\frac{1}{2}}-{\bf
w}_{i-\frac{1}{2},j-\frac{1}{2}}}{\tau}+B\frac{\widetilde{\mathcal
W}_{i,j-\frac{1}{2}}-\widetilde{\mathcal
W}_{i-1,j-\frac{1}{2}}}{h}+C\frac{\widetilde{\mathcal
W}_{i-\frac{1}{2},j}-\widetilde{\mathcal
W}_{i-\frac{1}{2},j-1}}{h}=\widetilde{f},
\end{equation}
with the initial condition $\mathbf{w}|_{t=0}=0$ where
$\mathbf{w}=\widehat{\upsilon}-\upsilon$, $\widetilde{\mathcal
W}=\widehat{\Upsilon}-\Upsilon$ and
 \begin{eqnarray*}
\widetilde{f}=(A^{(k)}-A)\frac{{\bf
\upsilon}^{i-\frac{1}{2},j-\frac{1}{2}}-{\bf
\upsilon}_{i-\frac{1}{2},j-\frac{1}{2}}}{\tau}+(B^{(k)}-B)\frac{
\Upsilon_{i,j-\frac{1}{2}}-
\Upsilon_{i-1,j-\frac{1}{2}}}{h}\\+(C^{(k)}-C)\frac{
{\Upsilon}_{i-\frac{1}{2},j}-
{\Upsilon}_{i-\frac{1}{2},j-1}}{h}.
 \end{eqnarray*}
By definition of $\widehat{\Upsilon}$ and $\Upsilon$, they are
constructed by means of Riemannian   invariants defined by the
``exact'' matrices $A,B,C$ and by their rational approximations
$A^{(k)},B^{(k)},C^{(k)}$, respectively. To work with the system~
\eqref{tilde} as with a difference scheme approximating the Cauchy
problem~\eqref{sist_1}, we would like to have in the left-hand part
of~\eqref{tilde} the ``large values''  ${\mathcal W}$ defined by
means of Riemannian invariants for $A,B,C$ and
$\widehat{\upsilon}-\upsilon$, rather then
$\widetilde{\mathcal{W}}$. Recall from Subsection~\ref{condif} (stage 2) that
the Riemannian invariants are constructed by means of eigenvectors
of the corresponding matrix pencils. Using the fact that the
cardinalities of their spectra and the spectrum of $A$ are known as
inputs, by \cite[Theorem~2]{zb04} we have that the eigenvectors
are computable, i.e.\ they can be chosen in such a way that, for an
absolute constant $c$,
 \begin{equation}\label{tmatr}
 ||T_{(A,B)}-T_{(A^{(k)},B^{(k)})}||_2\leq \frac{c}{2^k},\ ||T_{(A,C)}-T_{(A^{(k)},C^{(k)})}||_2\leq \frac{c}{2^k},
  \end{equation}
where by $T_{(A,B)}$ we denote the matrix $T_x$ defined in~\eqref{T}, and, in a
similar way for other pairs of matrices (cf.\ \cite[Theorem~5]{ss09}). The results of \cite{zb04} can be applied to this case since the procedure of finding $T_{(A,B)}$ consists of two spectral decompositions for symmetric matrices. Moreover, the fact that the matrix pencils $\lambda A-B$ and $\lambda A-C$ have no zero eigenvalues guarantees that the numbers of positive and negative eigenvalues in (13) are the same for $\lambda A-B$ and $\lambda A^{(k)}-B^{(k)}$, as well as for $\lambda A-C$ and $\lambda A^{(k)}-C^{(k)}$. Thus we can rewrite~\eqref{tilde} as
 \begin{eqnarray*}
 A\frac{{\bf w}^{i-\frac{1}{2},j-\frac{1}{2}}-{\bf
w}_{i-\frac{1}{2},j-\frac{1}{2}}}{\tau}+B\frac{{\mathcal
W}_{i,j-\frac{1}{2}}-{\mathcal
W}_{i-1,j-\frac{1}{2}}}{h}+C\frac{{\mathcal
W}_{i-\frac{1}{2},j}-{\mathcal W}_{i-\frac{1}{2},j-1}}{h}={f},
 \end{eqnarray*}
 where ${\mathcal W}$ are the desired ``large values'' for $w$ while $f$ differs from $\tilde{f}$ on a value involving the norms from~\eqref{tmatr} multiplied by difference derivatives analogous to the ones listed below in~\eqref{tmatr1}.

By the stability condition,
$||\widehat{\upsilon}-\upsilon||_{sL_2}\leq c||f||_{L_2}$ for some
constant $c$ depending only on $M_A, M_{\varphi}$ \cite{go71,go76,evans}. By
formal differentiation of the scheme (cf.\ \cite{go71}, similar
arguments in the proof of the existence theorem for one-dimensional symmetric hyperbolic systems in the canonical form) it is possible to check
that any of the norms
 \begin{equation}\label{tmatr1}
 ||\frac{{\bf
\upsilon}^{i-\frac{1}{2},j-\frac{1}{2}}-{\bf
\upsilon}_{i-\frac{1}{2},j-\frac{1}{2}}}{\tau}||_{sL_2},\
||\frac{\Upsilon_{i,j-\frac{1}{2}}-
\Upsilon_{i-1,j-\frac{1}{2}}}{h}||_{sL_2}, \ ||\frac{
\Upsilon_{i-\frac{1}{2},j}-
\Upsilon_{i-\frac{1}{2},j-1}}{h}||_{sL_2}
 \end{equation}
 is below the difference
derivatives of $\varphi_k$, which are below a constant depending only
on $M_A,M_{\varphi}$. We provide necessary arguments for the first norm, the arguments for the others being similar.

Direct computations (writing the system of difference equations~
\eqref{sist_dif_hat} for two neighbour time levels, subtracting
them and dividing by $\tau$) show that the grid function
 $$\frac{\Delta\upsilon^{(h)}}{\Delta t}=
 \left\{\frac{\upsilon^{l+1}_{i-\frac{1}{2},j-\frac{1}{2}}
 -\upsilon^l_{i-\frac{1}{2},j-\frac{1}{2}}}{\tau}\right\},
 $$
consisting of the ``difference derivatives'' of
$v^{(h)}$, meets the  system of difference
equations~\eqref{tmatr2} below.

For $t=l\tau$ we obtain:
 $$\begin{cases}
 A^{(k)}\frac{\upsilon^{l+1}_{i-\frac{1}{2},j-\frac{1}{2}}
-\upsilon^l_{i-\frac{1}{2},j-\frac{1}{2}}}{\tau}
+B^{(k)}\frac{\Upsilon^{l}_{i,j-\frac{1}{2}}
-\Upsilon^l_{i-1,j-\frac{1}{2}}}{h}+C^{(k)}\frac{\Upsilon^{l}_{i,j-\frac{1}{2}}
-\Upsilon^l_{i-1,j-\frac{1}{2}}}{h}=0,\\
\upsilon^0_{i-\frac{1}{2},j}=
(\varphi_k)_{i-\frac{1}{2},j-1}.
\end{cases}
 $$
For $t=(l+1)\tau$ we obtain
 $$\begin{cases} A^{(k)}\frac{\upsilon^{l+2}_{i-\frac{1}{2},j-\frac{1}{2}}
 -\upsilon^{l+1}_{i-\frac{1}{2},j-\frac{1}{2}}}{\tau}
+B^{(k)}\frac{\Upsilon^{l+1}_{i,j-\frac{1}{2}}
-\Upsilon^{l+1}_{i-1,j-\frac{1}{2}}}{h}+C^{(k)}\frac{\Upsilon^{l+1}_{i-\frac{1}{2},j}
-\Upsilon^{l+1}_{i-\frac{1}{2},j-1}}{h}=0,\\
\upsilon^1_{i-\frac{1}{2},j-\frac{1}{2}}=
(\varphi_k)_{i-\frac{1}{2},j-\frac{1}{2}}-(A^{(k)})^{-1}\tau \left(
B^{(k)}\frac{\Upsilon^0_{i,j-\frac{1}{2}}
-\Upsilon^0_{i-1,j-\frac{1}{2}}}{h}+C^{(k)}\frac{\Upsilon^0_{i-\frac{1}{2},j}
-\Upsilon^0_{i-\frac{1}{2},j-1}}{h}\right).
\end{cases}
 $$

 Subtracting the first system from the second one and dividing by $\tau$ we obtain:
 \begin{equation}\label{tmatr2}
 \begin{cases} A^{(k)}\frac{(\frac{\Delta \upsilon}{\Delta t})^{l+1}
-(\frac{\Delta \upsilon}{\Delta t})^{l}}{\tau}
+B^{(k)}\frac{(\frac{\Delta \Upsilon_x}{\Delta t})_{i}
-(\frac{\Delta \Upsilon_x}{\Delta
t})_{i-1}}{h}+C^{(k)}\frac{(\frac{\Delta \Upsilon_y}{\Delta
t})_{j}
-(\frac{\Delta \Upsilon_y}{\Delta t})_{j-1}}{h}=0,\\
(\frac{\Delta \upsilon}{\Delta t})^0= -(A^{(k)})^{-1} \left(
B^{(k)}\frac{\Upsilon^0_{i,j-\frac{1}{2}}-\Upsilon^0_{i-1,j-\frac{1}{2}}}{h}+C^{(k)}
\frac{\Upsilon^0_{i-\frac{1}{2},j}-\Upsilon^0_{i-\frac{1}{2},j-1}}{h}\right)
=:\psi_{i-\frac{1}{2}, j-\frac{1}{2}}
\end{cases}
 \end{equation}
 where
 \begin{eqnarray*}
 \left(\frac{\Delta \Upsilon_x}{\Delta t}\right)_i=\frac{
\Upsilon^{l+1}_{i,j-\frac{1}{2}}-
\Upsilon^{l}_{i,j-\frac{1}{2}}}{\tau},\quad
\left(\frac{\Delta \Upsilon_y}{\Delta t}\right)_j=(\frac{
\Upsilon^{l+1}_{i-\frac{1}{2},j}-
\Upsilon^{l}_{i-\frac{1}{2},j-1}}{\tau}),\\
\left(\frac{\Delta \upsilon}{\Delta t}\right)^l=\frac{
\upsilon^{l+1}_{i-\frac{1}{2},j-\frac{1}{2}}-
\upsilon^{l}_{i-\frac{1}{2},j-\frac{1}{2}}}{\tau}.
 \end{eqnarray*}
 Note that, since the operator of difference differentiation is linear, $\frac{\Delta \Upsilon_x}{\Delta t}$ and $\frac{\Delta \Upsilon_y}{\Delta t}$ are the ``large values'' for $\frac{\Delta \upsilon}{\Delta t}$.
The stability condition for this scheme looks as follows:
 $$
||\left\{\frac{\Delta \upsilon}{\Delta
t}\right\}_{i-\frac{1}{2}, j-\frac{1}{2}}||_{A^{(k)},sL_2}\leq
||\{\psi_{i-\frac{1}{2}, j-\frac{1}{2}}\}||_{A^{(k)},L_2}.
 $$

Since $A^{(k)},B^{(k)},C^{(k)}$ fast converge to $A,B,C$ respectively, it remains to estimate
 \begin{equation}\label{uh}
 ||\frac{\Delta U_x}{\Delta x}||:=||\left\{\frac{{\Upsilon}^0_{i, j-\frac{1}{2}}-
 {\Upsilon^0_{i-1,
 j-\frac{1}{2}}}}{h} \right\}||_{L_2},\quad ||\frac{\Delta U_y}{\Delta y}||:=
 ||\left\{\frac{{\Upsilon^0}_{i-\frac{1}{2},j}-
 {\Upsilon^0}_{i-\frac{1}{2},j-1}}{h} \right\}||_{L_2}.
 \end{equation}

Recall that the ``large values'' are calculated from the Riemannian
invariants of auxiliary one-dimensional systems, see Subsection
4.2: $\Upsilon^0_{i, j-\frac{1}{2}}=T_xV^0_i$, where
$V^0_i=v_{i\pm \frac{1}{2}}$, $V^0_{i-1}=v_{i-\frac{1}{2}}$ or
$v_{i-\frac{3}{2}}$, depending on the eigenvalues of the
corresponding matrix pencils (and in a similar way with
$\Upsilon^0_{i-\frac{1}{2},j}$). Note also that the eigenvectors were chosen to be orthonormal, hence $||T_x||_2=||T_x^{\ast}||_2=1$.
 Therefore,
  $$
||\frac{\Delta U_x}{\Delta x}||=||T_x^\ast\frac{\Delta V}{\Delta
x}||\leq 2||\frac{\Delta v}{\Delta x}||\leq c(M_A)||\frac{\Delta
\varphi_k}{\Delta x}||.
 $$
 The last estimate can be derived as an energy integral inequality for the auxiliary one-dimensional scheme like in~\eqref{ust} (see also
\cite[p.\ 78]{go76}, \cite{fri,gr,gv}).
 We have
   \begin{eqnarray*}
||\frac{\Delta \varphi_k}{\Delta x}||^2=\sum_{i,j}h^2 \left(
\frac{(\varphi_k)_{i+\frac{1}{2}, j-\frac{1}{2}}-
(\varphi_k)_{i-\frac{1}{2}, j-\frac{1}{2}}}{h}\right)^2=\\
\frac{1}{h^2}\sum_{i,j}h^2\langle (\varphi_k)_{i+\frac{1}{2},
j-\frac{1}{2}}- (\varphi_k)_{i-\frac{1}{2}, j-\frac{1}{2}},
(\varphi_k)_{i+\frac{1}{2}, j-\frac{1}{2}}-
(\varphi_k)_{i-\frac{1}{2}, j-\frac{1}{2}}\rangle.
 \end{eqnarray*}

 Adding and subtracting the expression $\varphi_{i+\frac{1}{2},
 j-\frac{1}{2}}- \varphi_{i-\frac{1}{2}, j-\frac{1}{2}}$ (where $\varphi$ is the ``exact'' initial function)
 within the scalar product, we obtain
 $$
||\frac{\Delta U_x}{\Delta x}||^2\leq\frac{2||\varphi_k-
\varphi||^2}{h^2}+ \sum_{i,j}h^2\langle
\frac{\varphi_{i+\frac{1}{2}, j-\frac{1}{2}}-
\varphi_{i-\frac{1}{2}, j-\frac{1}{2}}}{h},
\frac{\varphi_{i+\frac{1}{2}, j-\frac{1}{2}}-
\varphi_{i-\frac{1}{2}, j-\frac{1}{2}}}{h}\rangle.
 $$
 Since $\{\varphi_k\}$ fast converges to $\varphi$,
the first summand is below 2. Passing to the limit in the above
inequality when $h$ tends to 0 and taking into account that the
integration operator is computable we see that the second summand (which tends
to $\int_Q\langle \frac{\partial\varphi}{\partial
x},\frac{\partial\varphi}{\partial x}\rangle dxdy$) is uniformly
(on $h$) bounded by a universal constant depending only on
$M_\varphi$. Therefore, $||\frac{\Delta U_x}{\Delta x}||$ is
uniformly (on $h$) bounded by a universal constant depending only on $M_A$ and $M_\varphi$. The expression $||\frac{\Delta U_y}{\Delta y}||$ is estimated in a similar way.

 Since $||A^{(k)}-A||_2$, $||B^{(k)}-B||_2$,
and $||C^{(k)}-C||_2$ are below $\frac{1}{2^k}$, the desired
estimate for $||f||$ follows. This completes the proof of the estimate~\eqref{11}.

\section{Conclusion}\label{concl}

From the proof of Theorems~\ref{main},~\ref{main2} we see that the
same idea may be applied to a broader class of systems for  which
there is a stable difference scheme and which satisfy the existence
and uniqueness theorem. In particular, most probably, analogs of our results hold for systems~\eqref{sist_1},~\eqref{sist_2} with the variable coefficients $A,B_i$ depending on $t,x_j$. However,  precise specification of the whole
class of such systems and of corresponding difference schemes
remains an open question.

As noted in Section~\ref{state}, the wave equation~\eqref{wave} can be reduced
(not in a  unique way)  to a symmetric hyperbolic system, thus its
solution operator is computable which gives an affirmative answer to
the question of paper \cite{wz02}, for the case of  dissipative
boundary conditions and an initial function with uniformly bounded
derivatives. Note, however, that we consider not Sobolev $H^s$
spaces of generalized functions, but $C^k$ spaces of continuously
differentiable functions; the smoothness $k$, which is to be
assumed, depends on the particular problem under consideration. To
prove  computability of the generalized solutions it would be
probably suitable to use finite element methods.

 The restriction on the initial function seems to be
rather strong from Computable Analysis viewpoint, but it is very
natural for Numerical Analysis (though we have never seen it
explicitly formulated in  Numerical Analysis theorems).
Indeed, it is well known, that any initial (and right-hand part)
function can be represented in the form of a Fourier series
(consider for simplicity the one-dimensional case)
$\varphi(x)=\sum\limits_{n}a_ne^{inx}$, the differentiation of
which gives $\sum\limits_{n}na_ne^{inx}$, i.e.\ ``fast oscillating''
functions lead to large derivatives and hence large convergence
constants which can make the scheme not convergent on a real
computer.

A similar situation is with the additional assumptions of Theorem~\ref{main2} about the apriori knowledge  of the spectra of $A$ and of matrix pencils and the absence of zero eigenvalues: the assumptions correspond well to the experience of numerical analysts. Namely, violation of these assumptions may lead to computational instabilities. We currently do not know whether our results hold without these assumptions. Note, however, that the cardinalities are known (for physical reasons) for some important systems invariant under rotations \cite{gm98}.

Finally, we would like to point  out that it would be interesting
to study the computational complexity of the considered problems,
in the spirit of \cite{ko}.  The algorithms suggested in our
paper are very time- and space-consuming (in fact, our proofs provide no explicit complexity bounds at all, establishing only computability). Finding feasible algorithms to
compute the solution operators seems to be a  challenge.

\section*{Acknowledgement}
The work was partially supported  by
a Marie Curie
International Research Staff Exchange Scheme Fellowship within the
7th European Community Framework Programme and by RFBR grant 13-01-00015. 
The first author was partially supported by the Grant of Russian Federation for the State Support of Researches (Agreement No 14.B25.31.0029). The second author was partially supported by EU Research and Innovation Staff Exchange Programme ``Computing with Infinite Data''.

We are grateful to S.K.\ Godunov and  K.\ Weihrauch for
interesting discussions, to V.L.\ Miroshnichenko and O.V.\ 
Kudinov for bibliographical hints, and to the anonymous referees for the careful reading and valuable comments.

\end{document}